\documentclass[12pt]{article}
\usepackage{amstex,amssymb}
\begin{document}
\newtheorem{thm}{Theorem}[section]
\newtheorem{lem}[thm]{Lemma}
\newtheorem{prop}[thm]{Proposition}
\newtheorem{conj}[thm]{Conjecture}
\newtheorem{cor}[thm]{Corollary}
\newenvironment{dfn}{\medskip\refstepcounter{thm}
\noindent{\bf Definition \thesection.\arabic{thm}\ }}{\medskip}
\newenvironment{cond}{\medskip\refstepcounter{thm}
\noindent{\bf Condition \thesection.\arabic{thm}\ }}{\medskip}
\newenvironment{ex}{\medskip\refstepcounter{thm}
\noindent{\bf Example \thesection.\arabic{thm}\ }}{\medskip}
\newenvironment{proof}[1][,]{\medskip\ifcat,#1
\noindent{\it Proof.\ }\else\noindent{\it Proof of #1.\ }\fi}
{\hfill\rule{.5em}{.8em}\medskip}
\def\eq#1{{\rm(\ref{#1})}}
\def\Re{\mathop{\rm Re}}
\def\vol{\mathop{\rm vol}}
\def\Hol{{\textstyle\mathop{\rm Hol}}}
\def\id{\mathop{\rm id}}
\def\hcf{\mathop{\rm hcf}}
\def\Spin{\mathop{\rm Spin}}
\def\GL{\mathop{\rm GL}}
\def\SO{\mathop{\rm SO}}  
\def\SU{\mathop{\rm SU}}
\def\R{{\mathbb R}}
\def\Z{{\mathbb Z}}
\def\C{{\mathbb C}}
\def\CP{\mathbb{CP}}
\def\d{{\rm d}}
\def\sst{\scriptscriptstyle}
\def\w{\wedge}
\def\ra{\rightarrow}
\def\ha{{\textstyle{1\over2}}}
\def\Ahat{{\skew5\hat A}} 
\def\ms#1{\vert#1\vert^2}
\def\lnm#1#2{\Vert #1 \Vert_{L^{#2}}} 
\def\cnm#1#2{\Vert #1 \Vert_{C^{#2}}} 
\def\md#1{\vert #1 \vert}
\def\bcnm#1#2{\bigl\Vert #1 \bigr\Vert_{C^{#2}}} 
\def\bmd#1{\bigl\vert #1 \bigr\vert}
\def\an#1{\langle#1\rangle}
\title{A new construction of compact 8-manifolds with holonomy Spin(7)}
\author{Dominic Joyce}
\date{October 1999}
\maketitle

\section{Introduction}
\label{n1}

In Berger's classification \cite{Ber} of holonomy groups of
Riemannian manifolds there are two special cases, the exceptional
holonomy groups $G_2$ in 7 dimensions and $\Spin(7)$ in 8
dimensions. Bryant \cite{Bry} and Bryant and Salamon 
\cite{BrSa} showed that such metrics exist locally, and
wrote down explicit, complete metrics with holonomy $G_2$ 
and $\Spin(7)$ on noncompact manifolds. 

The first examples of metrics with holonomy $G_2$ and
$\Spin(7)$ on {\it compact} 7- and 8-manifolds were
constructed by the author in \cite{Joy1,Joy2,Joy3}. The 
survey paper \cite{Joy4} provides a good introduction 
to these constructions. Here is a brief description of the 
method used in \cite{Joy1} to construct compact 8-manifolds 
with holonomy $\Spin(7)$, divided into four steps.

\begin{itemize}
\setlength{\parsep}{0pt}
\setlength{\itemsep}{0pt}
\item[(a)] We start with a flat $\Spin(7)$-structure 
$(\Omega_0,g_0)$ on the 8-torus $T^8$, and a finite group 
$\Gamma$ of isometries of $T^8$ preserving $(\Omega_0,g_0)$. 
Then $T^8/\Gamma$ is an {\it orbifold}, a singular manifold 
with only quotient singularities.

\item[(b)] For certain $\Gamma$ one can resolve the singularities 
of $T^8/\Gamma$ in a natural way, using complex geometry. This 
gives a nonsingular, compact 8-manifold $M$, and a 
projection~$\pi:M\rightarrow T^8/\Gamma$.

\item[(c)] We write down a 1-parameter family of 
$\Spin(7)$-structures $(\Omega_t,g_t)$ on $M$ for $t\in(0,\epsilon)$,
such that $(\Omega_t,g_t)$ has small torsion when $t$ is small, and
converges to the singular $\Spin(7)$-structure $\pi^*(\Omega_0,g_0)$
as~$t\rightarrow 0$.

\item[(d)] We prove using analysis that for small $t$, 
the $\Spin(7)$-structure $(\Omega_t,g_t)$ can be deformed to a 
nearby $\Spin(7)$-structure $(\tilde\Omega_t,\tilde g_t)$ on $M$, 
with zero torsion. Then $\tilde g_t$ has holonomy~$\Spin(7)$. 
\end{itemize}

This paper describes a new method for constructing compact 
8-manifolds with holonomy $\Spin(7)$, in which one starts not 
with a torus $T^8$ but with a {\it Calabi--Yau $4$-orbifold} $Y$ 
with isolated singular points $p_1,\ldots,p_k$. We use algebraic 
geometry to find a number of suitable complex orbifolds $Y$, 
which in the simplest cases are hypersurfaces in {\it weighted 
projective spaces}~$\CP^5_{a_0,\ldots,a_5}$. 

Then, instead of a finite group $\Gamma$, we suppose we have an 
antiholomorphic, isometric involution $\sigma:Y\ra Y$, whose only 
fixed points are $p_1,\ldots,p_k$. This involution does not 
preserve the $\SU(4)$-structure on $Y$, but it does preserve 
the induced $\Spin(7)$-structure. We think of $\sigma$ as breaking 
the structure group of $Y$ from $\SU(4)$ down to $\Spin(7)$. Define 
$Z=Y/\an{\sigma}$. Then $Z$ is an orbifold with isolated singular 
points $p_1,\ldots,p_k$, and the Calabi--Yau structure on $Y$
induces a torsion-free $\Spin(7)$-structure on~$Z$.

If the singularities of $Z$ are of a suitable kind, we can
resolve them to get a compact 8-manifold $M$ with holonomy 
$\Spin(7)$, as in steps (b)--(d) above. To perform the
resolution we need to find {\it Asymptotically Locally 
Euclidean $\Spin(7)$-manifolds} corresponding to the 
singularities of $Z$, which are a special class of
noncompact $\Spin(7)$-manifolds asymptotic to quotient
singularities~$\R^8/G$.

Our construction then yields new examples of compact 8-manifolds 
$M$ with holonomy $\Spin(7)$. We calculate the Betti numbers 
$b^k(M)$ in each case. They turn out to be rather different to 
the Betti numbers arising from the previous construction in
\cite{Joy1}. In particular, in this new construction the 
middle Betti number $b^4$ tends to be rather large, as big as 
$11\,662$ in one example, whereas the manifolds of \cite{Joy1} 
all satisfied~$b^4\le 162$.

Sections \ref{n2} and \ref{n3} introduce the holonomy group
$\Spin(7)$ and Calabi--Yau orbifolds, and \S\ref{n4} defines the 
idea of ALE $\Spin(7)$-manifold, and gives a number of examples. 
Section \ref{n5} then proves our main result, that given a 
Calabi--Yau 4-orbifold $Y$ and an antiholomorphic involution
$\sigma:Y\ra Y$ satisfying certain conditions, we can construct
a compact 8-manifold $M$ with holonomy~$\Spin(7)$.

We explain in \S\ref{n6} how to use the construction in
practice, and ways of computing the Betti numbers of the
resulting 8-manifolds $M$. Sections \ref{n7}--\ref{n10}
apply the construction to generate new examples of compact
8-manifolds with holonomy $\Spin(7)$, and we finish in
\S\ref{n11} with a discussion of our results. 

The material in this paper will be discussed in the author's 
book \cite{Joy5}, which pays much attention to the exceptional 
holonomy groups, and also gives a more sophisticated version of 
the original construction \cite{Joy1} of compact 8-manifolds 
with holonomy~$\Spin(7)$.

\section{Background on the holonomy group $Spin(7)$}
\label{n2}

We now collect together some facts we will need about the
holonomy group $\Spin(7)$, taken from the books by Salamon 
\cite[Ch.~12]{Sal} and the author \cite[Ch.~10]{Joy5}.
First we define $\Spin(7)$ as a subgroup of~$\GL(8,\R)$.

\begin{dfn} Let $\R^8$ have coordinates $(x_1,\dots,x_8)$. 
Write $\d{\bf x}_{ijkl}$ for the 4-form 
$\d x_i\w\d x_j\w\d x_k\w\d x_l$ on $\R^8$. Define a 4-form 
$\Omega_0$ on $\R^8$ by
\begin{equation}
\begin{split}
\Omega_0=&\d{\bf x}_{1234}+\d{\bf x}_{1256}
+\d{\bf x}_{1278}+\d{\bf x}_{1357}-\d{\bf x}_{1368}\\
-&\d{\bf x}_{1458}-\d{\bf x}_{1467}-\d{\bf x}_{2358}
-\d{\bf x}_{2367}-\d{\bf x}_{2457}\\
+&\d{\bf x}_{2468}+\d{\bf x}_{3456}
+\d{\bf x}_{3478}+\d{\bf x}_{5678}.
\end{split}
\label{Om0eq}
\end{equation}
The subgroup of $\GL(8,\R)$ preserving $\Omega_0$ is the holonomy 
group $\Spin(7)$. It is a compact, connected, simply-connected, 
semisimple, 21-dimensional Lie group, which is isomorphic as 
a Lie group to the double cover of $\SO(7)$. This group also 
preserves the orientation on $\R^8$ and the Euclidean metric 
$g_0=\d x_1^2+\cdots+\d x_8^2$ on~$\R^8$. 

Let $M$ be an 8-manifold. For each $p\in M$, define ${\mathcal A}_pM$ 
to be the subset of 4-forms $\Omega\in\Lambda^4T_p^*M$ for which there 
exists an isomorphism between $T_pM$ and $\R^8$ identifying $\Omega$ 
and the 4-form $\Omega_0$ of \eq{Om0eq}. Let ${\mathcal A}M$ be the 
bundle with fibre ${\mathcal A}_pM$ at each $p\in M$. Then 
${\mathcal A}M$ is a subbundle of $\Lambda^4T^*M$ with fibre 
$\GL(8,\R)/\Spin(7)$. It is not a vector subbundle, and has 
codimension 27 in $\Lambda^4T^*M$. We say that a 4-form $\Omega$ on 
$M$ is {\it admissible} if $\Omega\vert_p\in{\mathcal A}_pM$ for 
each~$p\in M$. 
\label{s7def}
\end{dfn}

There is a 1-1 correspondence between $\Spin(7)$-structures $Q$ and 
admissible 4-forms $\Omega\in C^\infty({\mathcal A}M)$ on $M$. Each 
$\Spin(7)$-structure $Q$ induces a 4-form $\Omega$, a metric $g$ 
and an orientation on $M$, corresponding to $\Omega_0$, $g_0$ and
the orientation on~$\R^8$. 

\begin{dfn} Let $M$ be an 8-manifold, $\Omega$ an admissible 
4-form on $M$, and $g$ the associated metric. We shall 
abuse notation by referring to the pair $(\Omega,g)$ as a 
$\Spin(7)$-{\it structure} on $M$. Let $\nabla$ be 
the Levi-Civita connection of $g$. We call $\nabla\Omega$ the 
{\it torsion} of $(\Omega,g)$, and we say that $(\Omega,g)$ is 
{\it torsion-free} if $\nabla\Omega=0$. A triple $(M,\Omega,g)$ 
is called a $\Spin(7)$-{\it manifold} if $M$ is an 8-manifold, 
and $(\Omega,g)$ a torsion-free $\Spin(7)$-structure on~$M$.
\label{s7sdef}
\end{dfn}

Let $(\Omega,g)$ be a $\Spin(7)$-structure on an 8-manifold $M$.
Then $(\Omega,g)$ is torsion-free if and only if $\d\Omega=0$. If
$(\Omega,g)$ is torsion-free then $g$ is Ricci-flat, and $M$ is 
spin and has a constant positive spinor. If $M$ is compact 
and $\Hol(g)=\Spin(7)$ then the positive Dirac operator 
$D_+:C^\infty(S_+)\ra C^\infty(S_-)$ has kernel $\R$ and cokernel 0.
Thus $D_+$ has index 1. 

But the index of $D_+$ is the $\Ahat$-{\it genus} $\Ahat(M)$,
and is given by
\begin{equation}
24\Ahat(M)=-1+b^1(M)-b^2(M)+b^3(M)+b^4_+(M)-2b^4_-(M),
\label{Ahateq}
\end{equation}
where $b^k=b^k(M)$ are the Betti numbers of $M$. Thus a compact 
8-manifold $M$ with holonomy $\Spin(7)$ must satisfy 
$b^3+b^4_+=b^2+b^4_-+25$. As in \cite[Th.~C]{Joy1}, 
one can use this to show:

\begin{thm} Let\/ $(M,\Omega,g)$ be a compact\/ $\Spin(7)$-manifold.
Then $\Hol(g)=\Spin(7)$ if and only if\/ $M$ is simply-connected,
and\/~$b^3+b^4_+=b^2+b^4_-+25$. 
\label{n2thm1}
\end{thm}

The following result \cite[Th.~D]{Joy1} describes the
moduli space of holonomy $\Spin(7)$ metrics.

\begin{thm} Let\/ $M$ be a compact\/ $8$-manifold admitting
metrics with holonomy $\Spin(7)$. Then the moduli space of
metrics with holonomy $\Spin(7)$ on $M$, up to diffeomorphisms
isotopic to the identity, is a smooth manifold of 
dimension~$1+b^4_-(M)$.
\label{n2thm2}
\end{thm}

Our next proposition follows from the ideas of~\cite[\S 10.6]{Joy5}.

\begin{prop} Let\/ $M$ be an $8$-manifold. Then there exists a 
tubular open neighbourhood ${\mathcal T}M$ of ${\mathcal A}M$ in 
$\Lambda^4T^*M$ which is a fibration over $M$, a smooth map of 
fibre bundles $\Theta:{\mathcal T}M\ra{\mathcal A}M$, and positive 
constants $\rho,C$, such that
\begin{itemize}
\item[{\rm(i)}] If\/ $(\Omega,g)$ is a $\Spin(7)$-structure and
$\xi$ a $4$-form on $M$ with\/ $\md{\xi-\Omega}_g\le\rho$,
then~$\xi\in C^\infty({\mathcal T}M)$.
\item[{\rm(ii)}] Suppose $(\Omega,g)$ is a $\Spin(7)$-structure on 
$M$, and\/ $\xi$ a $4$-form on $M$ with\/ $\md{\xi-\Omega}_g\le\rho$. 
Write $\Omega'=\Theta(\xi)$, and let\/ $(\Omega',g')$ be the associated 
$\Spin(7)$-structure. Then $\md{\xi-\Omega'}_{g'}\le\md{\xi-\Omega}_g$.
If\/ $(\Omega,g)$ is also torsion-free, 
then~$\bmd{\nabla'(\xi-\Omega')}_{g'}\le C\bmd{\nabla(\xi-\Omega)}_g$.
\end{itemize}
Here $\nabla,\nabla'$ are the Levi-Civita connections of\/
$g$ and\/ $g'$, and\/ $\md{\,.\,}_g$, $\md{\,.\,}_{g'}$ 
the norms defined using $g$ and~$g'$.
\label{n2prop}
\end{prop}

This is an entirely local result, involving calculations at a
point, and $\rho,C$ are independent of $M$. The inequality 
$\md{\xi-\Omega'}_{g'}\le\md{\xi-\Omega}_g$ in part (ii) 
should be understood as saying that $\Omega'=\Theta(\xi)$ is 
the $\Spin(7)$-form closest to $\xi$. That is, ${\mathcal T}M$ is 
a small open neighbourhood of ${\mathcal A}M$ in $\Lambda^4T^*M$, and
$\Theta$ is the projection from ${\mathcal T}M$ to the nearest point
in ${\mathcal A}M$. But as we have not fixed a metric on $M$, we do 
not have a way to measure distance in $\Lambda^4T^*M$, and so we 
use the metrics $g$, $g'$ associated to the $\Spin(7)$-forms
$\Omega,\Omega'$ to do this.

Our final result is proved in \cite[Th.~A \& Th.~B]{Joy1},
and also in~\cite[Ch.~13]{Joy5}.

\begin{thm} Let\/ $\lambda,\mu,\nu$ be positive constants. Then 
there exist positive constants $\kappa,K$ such that whenever
$0<t\le\kappa$, the following is true.

Let\/ $M$ be a compact\/ $8$-manifold, and\/ $(\Omega,g)$ a 
$\Spin(7)$-structure on $M$. Suppose that\/ $\phi$ is a 
smooth\/ $4$-form on $M$ with\/ $\d\Omega+\d\phi=0$, and 
\begin{itemize}
\setlength{\itemsep}{0pt}
\setlength{\parsep}{0pt}
\item[{\rm(i)}] $\lnm{\phi}2\le\lambda t^{9/2}$ 
and\/~$\lnm{\d\phi}{10}\le\lambda t$,
\item[{\rm(ii)}] the injectivity radius $\delta(g)$ 
satisfies $\delta(g)\ge\mu t$, and
\item[{\rm(iii)}] the Riemann curvature $R(g)$ 
satisfies~$\bcnm{R(g)}0\le\nu t^{-2}$.
\end{itemize}
Then there exists a smooth, torsion-free $\Spin(7)$-structure 
$(\tilde\Omega,\tilde g)$ on $M$ with\/
$\cnm{\tilde\Omega-\Omega}0\le Kt^{1/2}$.
\label{n2thm3}
\end{thm}

Here is how to interpret this result. As $\nabla\Omega=0$ if 
and only if $\d\Omega=0$ and $\d\phi+\d\Omega=0$, the torsion 
$\nabla\Omega$ is determined by $\d\phi$. Thus we can think of 
$\phi$ as a {\it first integral of the torsion} of $(\Omega,g)$. 
So $\lnm{\phi}2$ and $\lnm{\d\phi}{10}$ are both measures of the 
torsion of $(\Omega,g)$. As $t$ is small, part (i) of the theorem 
says that $(\Omega,g)$ has {\it small torsion} in a certain sense.

Parts (ii) and (iii) say that the injectivity radius of $g$
should not be too small, and its curvature not too large.
When a metric becomes singular, in general its injectivity
radius goes to zero and its curvature becomes infinite. So
we can interpret (ii) and (iii) as saying that $g$ is not
too close to being singular.

Thus, the theorem as a whole says that if the torsion of
$(\Omega,g)$ is small enough, and $g$ is not too singular,
then we can deform $(\Omega,g)$ to a nearby, torsion-free
$\Spin(7)$-structure $(\tilde\Omega,\tilde g)$ on $M$. We can
then use Theorem \ref{n2thm1} to show that if $M$ is
simply-connected and $b^3+b^4_+=b^2+b^4_-+25$ then
$\tilde g$ has holonomy~$\Spin(7)$.

We prove Theorem \ref{n2thm3} using analysis: we write
the condition that $(\tilde\Omega,\tilde g)$ be torsion-free 
as a nonlinear elliptic p.d.e., which can be approximated
by a linear elliptic p.d.e.\ when $\tilde\Omega-\Omega$ is 
small. Then we use tools such as Sobolev spaces, the Sobolev 
Embedding Theorem and elliptic regularity to show that
this nonlinear elliptic p.d.e.\ has a smooth solution.

\section{Calabi--Yau manifolds and orbifolds}
\label{n3}

We now give a brief introduction to Calabi--Yau geometry, and the 
relation between Calabi--Yau 4-folds and $\Spin(7)$-manifolds.
Some suitable references are Salamon \cite[Ch.~8]{Sal} and the
author~\cite[Ch.~6]{Joy5}.

\begin{dfn} A {\it Calabi--Yau manifold} or {\it orbifold} is 
a compact K\"ahler manifold or orbifold $(Y,J,g)$ of dimension $m$, 
with $\Hol(g)=\SU(m)$.
\label{cyodef}
\end{dfn}

Now Calabi--Yau manifolds and orbifolds are nearly the same thing 
as Ricci-flat K\"ahler manifolds and orbifolds, as we see in the next 
proposition. It follows from elementary properties of holonomy 
groups and K\"ahler geometry.

\begin{prop} Any Calabi--Yau orbifold $(Y,J,g)$ is Ricci-flat. 
Conversely, let\/ $(Y,J,g)$ be a compact Ricci-flat K\"ahler orbifold 
of dimension $m$, with singular set\/ $S$. Suppose that\/ $Y\setminus S$ 
is simply-connected and\/ $h^{p,0}(Y)=0$ for $0<p<m$. Then 
$\Hol(g)=\SU(m)$, so $Y$ is a Calabi--Yau orbifold.
\label{n3prop1}
\end{prop}

But using Yau's proof of the Calabi conjecture \cite{Yau}, one can 
show that suitable complex orbifolds admit Ricci-flat K\"ahler metrics.

\begin{thm} Let\/ $(Y,J)$ be a compact complex orbifold
admitting K\"ahler metrics, with\/ $c_1(Y)=0$. Then there is 
a unique Ricci-flat K\"ahler metric in each K\"ahler class on~$Y$. 
\label{cyrfthm}
\end{thm}

Now the action of $\SU(m)$ on $\C^m$ fixes the complex $m$-form
$\d z_1\w\cdots\w\d z_m$. It follows by general principles of
Riemannian holonomy that any Riemannian manifold or orbifold with
holonomy $\SU(m)$ admits a complex $m$-form $\theta$ corresponding to 
$\d z_1\w\cdots\w\d z_m$ which is constant under the Levi-Civita
connection $\nabla$. So we get:

\begin{prop} Let\/ $(Y,J,g)$ be a Calabi--Yau manifold or orbifold of 
dimension $m$, with K\"ahler form $\omega$. Then there exists a constant\/ 
$(m,0)$-form $\theta$ on $Y$, such that near every point\/ $p\in Y$ we 
can choose complex coordinates $(z_1,\ldots,z_m)$ in which
\begin{equation}
\begin{split}
\!\!\!\!\!\!\!\!\!\!\!\!\!\!\!\!
g=\ms{\d z_1}+\cdots+\ms{\d z_m},\quad
\omega&={i\over 2}(\d z_1\w\d\bar z_1+\cdots+\d z_m\w\d\bar z_m),\\
\text{and}\quad\theta&=\d z_1\w\cdots\w\d z_m
\end{split}
\label{sumeq}
\end{equation}
at\/ $p$. This form $\theta$ is unique up to multiplication by
$e^{i\phi}$ for some~$\phi\in[0,2\pi)$.
\label{n3prop2}
\end{prop}

We call $\theta$ the {\it holomorphic volume form} of $Y$. Now we 
restrict our attention to complex dimension 4. Here is a
criterion for a complex 4-orbifold to be Calabi--Yau.

\begin{prop} Let\/ $(Y,J)$ be a compact complex $4$-orbifold
with\/ $c_1(Y)=0$, admitting K\"ahler metrics. Suppose $Y\setminus S$ 
is simply-connected, where $S$ is the singular set of\/ $Y$, and\/ 
$h^{2,0}(Y)=0$. Then each K\"ahler class on $Y$ contains a unique 
metric $g$ such that\/ $(Y,J,g)$ is a Calabi--Yau $4$-orbifold.
\label{n3prop3}
\end{prop}

\begin{proof} As $\pi_1(Y\setminus S)=0$ we have $b^1(Y)=0$, so that 
$h^{1,0}(Y)=0$. Since $\pi_1(Y\setminus S)=0$ and $c_1(Y)=0$ the 
canonical bundle $K_Y$ of $Y$ is trivial, and this implies
that $h^{p,0}(Y)=h^{4-p,0}(Y)$. Thus $h^{3,0}(Y)=0$. But we
are given that $h^{2,0}(Y)=0$. Hence $h^{p,0}(Y)=0$ for
$0<p<4$, and the proposition follows from Proposition 
\ref{n3prop1} and Theorem~\ref{cyrfthm}.
\end{proof}

A Calabi--Yau 4-fold $Y$ has holonomy $\SU(4)$, and so carries 
a natural torsion-free $\SU(4)$-structure. Since $\SU(4)\subset
\Spin(7)\subset\SO(8)$, this $\SU(4)$-structure induces a 
$\Spin(7)$-structure on $Y$, which is also torsion-free.

\begin{prop} Suppose $(Y,J,g)$ is a Calabi--Yau $4$-orbifold, 
with K\"ahler form $\omega$ and holomorphic volume form $\theta$. Define 
a $4$-form $\Omega$ on $Y$ by $\Omega=\ha\omega\w\omega+\Re(\theta)$. 
Then $(\Omega,g)$ is a torsion-free $\Spin(7)$-structure on~$Y$.
\label{n3prop4}
\end{prop}

\begin{proof} Let $p$ be a point in $Y$. Then by Proposition
\ref{n3prop2} we can choose complex coordinates $(z_1,\ldots,z_4)$
near $p$ such that $g,\omega$ and $\theta$ are given by \eq{sumeq} at $p$,
with $m=4$. Define real coordinates $(x_1,\ldots,x_8)$ on $Y$ near
$p$ such that $(z_1,\ldots,z_4)=(x_1+ix_2,x_3+ix_4,x_5+ix_6,x_7+ix_8)$. 
Then from \eq{sumeq} we see that $g$, $\omega$ and $\Re(\theta)$ are given
at $p$ by 
\begin{gather*}
g=\d x_1^2+\cdots+\d x_8^2,\quad
\omega=\d{\bf x}_{12}+\d{\bf x}_{34}+\d{\bf x}_{56}+\d{\bf x}_{78}
\quad\text{and}\\
\Re(\theta)=\d{\bf x}_{1357}\!-\d{\bf x}_{1368}\!-\d{\bf x}_{1458}\!
-\d{\bf x}_{1467}\!-\d{\bf x}_{2358}\!-\d{\bf x}_{2367}\!
-\d{\bf x}_{2457}\!+\d{\bf x}_{2468},
\end{gather*}
where $\d{\bf x}_{ij\ldots l}=\d x_i\w\d x_j\w\cdots\w\d x_l$. 

It follows from this equation that $\Omega=\ha\omega\w\omega+\Re(\theta)$
coincides with the 4-form $\Omega_0$ defined in \eq{Om0eq}. As
this holds for all $p\in Y$, we see that $(\Omega,g)$ is a
$\Spin(7)$-{\it structure} on $Y$, in the sense of Definition
\ref{s7sdef}. Now $\nabla\omega=\nabla\theta=0$, where $\nabla$ is
the Levi-Civita connection of $g$, and so $\nabla\Omega=0$. But
$\nabla\Omega$ is the {\it torsion} of $(\Omega,g)$, so that $(\Omega,g)$
is torsion-free, as we want.
\end{proof}

Thus Calabi-Yau 4-folds are also $\Spin(7)$-manifolds.

\section{ALE $\Spin(7)$-manifolds}
\label{n4}

{\it ALE manifolds}, or {\it Asymptotically Locally Euclidean 
manifolds}, are a class of noncompact Riemannian manifolds with one end
modelled asymptotically on a quotient singularity~$\R^n/G$.

\begin{dfn} Let $G$ be a finite subgroup of $\SO(n)$ which acts 
freely on $\R^n\setminus\{0\}$. Let $X$ be a noncompact $n$-manifold
and $\pi:X\ra\R^n/G$ a continuous, surjective map, such that
$\pi^{-1}(0)$ is a compact subset of $X$, and $\pi:X\setminus\pi^{-1}(0)
\ra(\R^n/G)\setminus\{0\}$ is a diffeomorphism. Then we call $(X,\pi)$ 
a {\it real resolution} of~$\R^n/G$.

A metric $g$ on $X$ is called {\it Asymptotically Locally 
Euclidean}, or {\it ALE}, if 
\begin{equation*}
\nabla^l(\pi_*(g)-g_0)=O(r^{-n-l})
\quad\text{on $\{x\in\R^n/G:r(x)>R\}$, for all $l\ge 0$.}
\end{equation*}
Here $g_0$ is the Euclidean metric on $\R^n/G$, $r$ is the radius 
function on $\R^8/G$, and $R>0$ is a constant. We say that $(X,g)$
is {\it asymptotic to}~$\R^n/G$.
\label{aledef}
\end{dfn}

One reason ALE manifolds are interesting is that if you have
an ALE manifold $(X,g_{\sst X})$ asymptotic to $\R^n/G$, and
a compact Riemannian orbifold $(Y,g_{\sst Y})$ with isolated 
singularities modelled on $\R^n/G$, then you can glue $X$ and 
$Y$ together to get a nonsingular, compact Riemannian manifold 
$(M,g_{\sst M})$. We think of this as resolving the 
singularities of $Y$ using~$X$.

This technique is particularly valuable when $X$ and $Y$ 
both have special holonomy, so that $\Hol(g_{\sst X})$ and
$\Hol(g_{\sst Y})$ both lie in some holonomy group 
$H\subset\SO(n)$, as then we can hope to construct a metric
$g_{\sst M}$ on $M$ with $\Hol(g_{\sst M})\subseteq H$. So ALE 
manifolds $(X,g_{\sst X})$ with $\Hol(g_{\sst X})\subseteq H$
are ingredients in a construction for compact manifolds with
holonomy~$H$.

In fact the only interesting candidates for the holonomy group $H$ 
are $U(m)$ and $\SU(m)$ for $m\ge 2$, and $\Spin(7)$. Kronheimer 
\cite{Kro1,Kro2} constructed and classified all ALE 4-manifolds 
with holonomy $\SU(2)$. Calabi \cite[p.~285]{Cal} found an 
explicit family of ALE manifolds with holonomy $\SU(m)$ asymptotic 
to $\C^m/\Z_m$, and more generally the author \cite{Joy6}, 
\cite[Ch.~8]{Joy5} gave existence theorems for ALE manifolds 
with holonomy $\SU(m)$. No examples of ALE 8-manifolds with 
holonomy $\Spin(7)$ are known, at the time of writing.

However, we can construct compact 8-manifolds with holonomy
$\Spin(7)$ using only ALE 8-manifolds whose holonomy is a proper 
subgroup of $\Spin(7)$ such as $\SU(4)$ or $\Z_2\ltimes\SU(4)$,
and many examples of these can be found using the results of
\cite{Joy6}. To discuss these, it is useful to define the idea
of {\it ALE\/ $\Spin(7)$-manifold}, as in~\cite[Ch.~13]{Joy5}.

\begin{dfn} Let $G$ be a finite subgroup of $\Spin(7)$ which acts 
freely on $\R^8\setminus\{0\}$, let $(X,\pi)$ be a real resolution of
$\R^8/G$, and $(\Omega,g)$ a torsion-free $\Spin(7)$-structure on $X$. 
We call $(X,\Omega,g)$ an {\it ALE\/ $\Spin(7)$-manifold} if 
\begin{equation*}
\nabla^l(\pi_*(\Omega)-\Omega_0)=O(r^{-8-l})
\quad\text{on $\{x\in\R^8/G:r(x)>R\}$, for all $l\ge 0$.}
\end{equation*}
Here $\Omega_0$ is the $\Spin(7)$ 4-form on $\R^8/G$ given in 
\eq{Om0eq}, $r$ the radius function on $\R^8/G$, and $R>0$ a constant.
\label{ales7def}
\end{dfn}

In the rest of the section we give some examples of ALE 
$\Spin(7)$-manifolds.

\subsection{An example of an ALE $\Spin(7)$-manifold}
\label{n41}

We define a finite group $G\subset\Spin(7)$, such that $\R^8/G$ 
has an isolated singularity at 0, and construct two topologically 
distinct ALE $\Spin(7)$-manifolds $(X_1,\Omega_1,g_1)$ and 
$(X_2,\Omega_2,g_2)$ asymptotic to $\R^8/G$. These will be used 
in \S\ref{n5} as part of a construction of compact 8-manifolds 
with holonomy~$\Spin(7)$.

Let $\R^8$ have coordinates $(x_1,\ldots,x_8)$ and 
$\Spin(7)$-structure $(\Omega_0,g_0)$, as in Definition 
\ref{s7def}. Use the complex coordinates
\begin{equation*}
(z_1,z_2,z_3,z_4)=(x_1+ix_2,x_3+ix_4,x_5+ix_6,x_7+ix_8)
\end{equation*}
to identify $\R^8$ with $\C^4$. Then $g_0=\ms{\d z_1}+\cdots+\ms{\d z_4}$, 
and $\Omega_0=\ha\omega_0\w\omega_0+\Re(\theta_0)$, where $\omega_0$ 
is the K\"ahler form of $g_0$ and $\theta_0=\d z_1\w\cdots\w\d z_4$ 
the complex volume form on~$\C^4$. 

Define $\alpha,\beta:\C^4\ra\C^4$ by
\begin{equation}
\begin{split}
\!\!\!\!\!\!\!\!\!\!\!\!\!\!\!\!\!\!
\alpha:(z_1,\ldots,z_4)&\mapsto(iz_1,iz_2,iz_3,iz_4),\\
\!\!\!\!\!\!\!\!\!\!\!\!\!\!\!\!\!\!
\beta:(z_1,\ldots,z_4)&\mapsto(\bar z_2,-\bar z_1,\bar z_4,-\bar z_3).
\end{split}
\label{n4abeq1}
\end{equation}
Then $\alpha\in\SU(4)\subset\Spin(7)$ and $\beta\in\Spin(7)$, and 
$\alpha,\beta$ satisfy $\alpha^4=\beta^4=1$, $\alpha^2=\beta^2$ and 
$\alpha\beta=\beta\alpha^3$. Let $G=\an{\alpha,\beta}$. Then $G$ is 
a finite nonabelian subgroup of $\Spin(7)$ of order 8 which acts 
freely on~$\R^8\setminus\{0\}$. 

Now $\C^4/\an{\alpha}$ is a complex singularity, as $\alpha\in\SU(4)$.
Let $(Y_1,\pi_1)$ be the blow-up of $\C^4/\an{\alpha}$ at 0. Then 
$Y_1$ is the unique crepant resolution of $\C^4/\an{\alpha}$. The 
action of $\beta$ on $\C^4/\an{\alpha}$ lifts to a free antiholomorphic 
map $\beta:Y_1\ra Y_1$ with $\beta^2=1$. Define $X_1=Y_1/\an{\beta}$. 
Then $X_1$ is a nonsingular 8-manifold, and the projection 
$\pi_1:Y_1\ra\C^4/\an{\alpha}$ pushes down to~$\pi_1:X_1\ra\R^8/G$. 

By \cite[Th.~3.3, Th.~3.4]{Joy6} there exist ALE K\"ahler metrics 
$g_1$ on $Y_1$ with holonomy $\SU(4)$, which were in fact written down 
explicitly by Calabi \cite[p.~285]{Cal}. Each such $g_1$ is invariant 
under the action of $\beta$ on $Y_1$. Let $\omega_1$ be the K\"ahler 
form of $g_1$, and $\theta_1=\pi_1^*(\theta_0)$ the holomorphic volume 
form on $Y_1$. Then Proposition \ref{n3prop4} defines a torsion-free 
$\Spin(7)$-structure $(\Omega_1,g_1)$ on $Y_1$ 
with~$\Omega_1=\ha\omega_1\w\omega_1+\Re(\theta_1)$. 

As $\beta^*(\omega_1)=-\omega_1$ and $\beta^*(\theta_1)=\bar\theta_1$, we 
see that $\beta$ preserves $(\Omega_1,g_1)$. Thus $(\Omega_1,g_1)$ pushes 
down to a torsion-free $\Spin(7)$-structure $(\Omega_1,g_1)$ on 
$X_1$. Then $(X_1,\Omega_1,g_1)$ is an {\it ALE\/ $\Spin(7)$-manifold} 
asymptotic to $\R^8/G$. The Betti numbers of $X_1$ are 
$b^1=b^2=b^3=0$ and $b^4=1$, and~$\pi_1(X_1)=\Z_2$.

\subsection{A second ALE $\Spin(7)$-manifold asymptotic to $\R^8/G$}
\label{n42}

Define new complex coordinates $(w_1,\ldots,w_4)$ on $\R^8$ by
\begin{equation*}
(w_1,w_2,w_3,w_4)=(-x_1+ix_3,x_2+ix_4,-x_5+ix_7,x_6+ix_8).
\end{equation*}
Then $g_0=\ms{\d w_1}+\cdots+\ms{\d w_4}$ and 
$\Omega_0=\ha\omega_0'\w\omega_0'+\Re(\theta_0')$, where $\omega_0'$ is 
the K\"ahler form of $g_0$ w.r.t.\ the complex structure induced by the 
$w_j$, and $\theta_0'=\d w_1\w\cdots\w\d w_4$ is the complex volume 
form on~$\C^4$. 

As the action of $\SU(4)$ on $\R^8=\C^4$ induced by the $w_j$ 
preserves $g_0,\omega_0'$ and $\theta_0'$, it preserves $(\Omega_0,g_0)$. 
Thus the action of $\SU(4)$ on $\R^8$ compatible with the 
coordinates $w_j$ is a subgroup of $\Spin(7)$. Note that this 
is a {\it different}\/ $\SU(4)$ subgroup of $\Spin(7)$ to that 
considered above, induced by the $z_j$. In the coordinates $w_j$, 
we find that $\alpha,\beta$ act by
\begin{equation}
\begin{split}
\alpha:(w_1,\ldots,w_4)&\mapsto(\bar w_2,-\bar w_1,\bar w_4,-\bar w_3),\\
\beta:(w_1,\ldots,w_4)&\mapsto(iw_1,iw_2,iw_3,iw_4).
\end{split}
\label{n4abeq2}
\end{equation}
Observe that \eq{n4abeq1} and \eq{n4abeq2} are the same, except
that the r\^oles of $\alpha,\beta$ are reversed. Therefore we can use 
the ideas above again.

Let $Y_2$ be the crepant resolution of $\C^4/\an{\beta}$. The action 
of $\alpha$ on $\C^4/\an{\beta}$ lifts to a free antiholomorphic
involution of $Y_2$. Let $X_2=Y_2/\an{\alpha}$. Then $X_2$ is nonsingular, 
and as above there exists a torsion-free $\Spin(7)$-structure 
$(\Omega_2,g_2)$ on $X_2$, making $(X_2,\Omega_2,g_2)$ into an ALE 
$\Spin(7)$-manifold asymptotic to~$\R^8/G$.

Now $(X_1,\Omega_1,g_1)$, $(X_2,\Omega_2,g_2)$ are clearly isomorphic 
as $\Spin(7)$-manifolds, but they should be regarded as 
{\it topologically distinct} ALE manifolds, because the 
isomorphism between them acts nontrivially on $\R^8/G$. 
Thus, we have found two topologically distinct ALE 
$\Spin(7)$-manifolds $(X_1,\Omega_1,g_1)$, $(X_2,\Omega_2,g_2)$ 
asymptotic to the same singularity~$\R^8/G$.

\subsection{Other examples of ALE $\Spin(7)$-manifolds}
\label{n43}

We can use the ideas above to construct other ALE 
$\Spin(7)$-manifolds too. Here we very briefly describe
two infinite families of ALE $\Spin(7)$-manifolds $X_1^n$, 
$X_2^n$ for $n=1,3,5,\ldots$. For simplicity they will not be 
used in the rest of the paper, although they easily could be.

Identify $\R^8$ and $\C^4$ as in \S\ref{n41}. Let $n\ge 1$ 
be an odd integer, and define $\alpha,\beta,\gamma:\C^4\ra\C^4$ by
\begin{align*}
\alpha:(z_1,\ldots,z_4)&\mapsto({\rm e}^{2\pi i/n}z_1,
{\rm e}^{-2\pi i/n}z_2,{\rm e}^{2\pi i/n}z_3,{\rm e}^{-2\pi i/n}z_4), \\
\beta:(z_1,\ldots,z_4)&\mapsto(iz_1,iz_2,iz_3,iz_4),\\
\gamma:(z_1,\ldots,z_4)&\mapsto(\bar z_2,-\bar z_1,\bar z_4,-\bar z_3).
\end{align*}
Then $\alpha,\beta\in\SU(4)$ and $\gamma\in\Spin(7)$, and 
$G^n=\an{\alpha,\beta,\gamma}$ is a finite nonabelian subgroup of 
$\Spin(7)$ of order $8n$ which acts freely on $\R^8\setminus\{0\}$. 
Note that $G^1$ coincides with the group $G$ of \S\ref{n41}-\S\ref{n42}.

We can construct a family of ALE $\Spin(7)$-manifolds asymptotic to
$\R^8/G^n$ as follows. The complex singularity $\C^4/\an{\alpha,\beta}$
has a unique crepant resolution $Y_1^n$, which can be described explicitly
using toric geometry. The action of $\gamma$ on $\C^4/\an{\alpha,\beta}$ 
lifts to a free antiholomorphic involution $\gamma:Y_1^n\ra Y_1^n$, so 
that $X_1^n=Y_1^n/\an{\gamma}$ is a nonsingular 8-manifold with a 
projection~$\pi_1^n:X_1^n\ra\R^8/G^n$.

By the results of \cite{Joy6}, there exist ALE K\"ahler metrics 
$g_1^n$ on $Y_1^n$ with holonomy $\SU(4)$. We can choose $g_1^n$ to 
be $\gamma$-invariant, and then the induced $\Spin(7)$-structure 
$(\Omega_1^n,g_1^n)$ on $Y_1^n$ is also $\gamma$-invariant, and 
pushes down to $X_1^n$, making $(X_1^n,\Omega_1^n,g_1^n)$ into an 
ALE $\Spin(7)$-manifold asymptotic to $\R^8/G^n$. Using the ideas of 
\S\ref{n42}, we can also construct a second ALE $\Spin(7)$-manifold
$(X_2^n,\Omega_2^n,g_2^n)$ asymptotic to~$\R^8/G^n$.

\section{Proof of the construction}
\label{n5}

Starting with a Calabi--Yau 4-orbifold $Y$ with isolated
singularities of a certain kind, and an antiholomorphic involution 
$\sigma$ on $Y$, we will now construct a compact 8-manifold $M$ by 
resolving $Z=Y/\an{\sigma}$, and prove that there exist torsion-free 
$\Spin(7)$-structures $(\tilde\Omega,\tilde g)$ on $M$, which have 
holonomy $\Spin(7)$ if $M$ is simply-connected. 

\subsection{A class of $\Spin(7)$-orbifolds $Z$} 
\label{n51}

We set out below the ingredients in our construction, and the
assumptions they must satisfy.

\begin{cond} Let $(Y,J)$ be a compact complex 4-orbifold with 
$c_1(Y)=0$, admitting K\"ahler metrics. Let $\sigma$ be an antiholomorphic 
involution on $Y$. That is, $\sigma:Y\ra Y$ is a diffeomorphism satisfying 
$\sigma^2=\id$ and $\sigma^*(J)=-J$. Define $\alpha:\C^4\ra\C^4$ by
\begin{equation}
\alpha:(z_1,z_2,z_3,z_4)\longmapsto(iz_1,iz_2,iz_3,iz_4).
\label{n5aleq}
\end{equation}
Then $\alpha^4=1$, so that $\an{\alpha}\cong\Z_4$, and $\C^4/\an{\alpha}$
has an isolated singular point at 0. We require that the singular 
set of $Y$ should be $k$ isolated points $p_1,\ldots,p_k$ for some 
$k\ge 1$, each modelled on $\C^4/\an{\alpha}$, and that the fixed set 
of $\sigma$ in $Y$ is exactly $\{p_1,\dots,p_k\}$. We also suppose that 
$Y\setminus\{p_1,\ldots,p_k\}$ is simply-connected, and~$h^{2,0}(Y)=0$.
\label{n5cond}
\end{cond}

In the rest of the section we assume that Condition \ref{n5cond} holds.

\begin{prop} There is a $\sigma$-invariant metric $g_{\sst Y}$ 
on $Y$ making $(Y,J,g_{\sst Y})$ into a Calabi--Yau orbifold.
We can choose the holomorphic volume form $\theta_{\sst Y}$ on 
$(Y,J,g_{\sst Y})$ such that\/ $\sigma^*(\theta_{\sst Y})=\bar\theta_{\sst Y}$.
Let\/ $(\Omega_{\sst Y},g_{\sst Y})$ be the torsion-free
$\Spin(7)$-structure on $Y$ from Proposition \ref{n3prop4}.
Then $(\Omega_{\sst Y},g_{\sst Y})$ is $\sigma$-invariant.
\label{n5prop1}
\end{prop}

\begin{proof} Let $g'$ be a K\"ahler metric on $Y$. Then $\sigma^*(g')$ 
is also a K\"ahler metric on $Y$, and so $g''=g'+\sigma^*(g')$ is a 
$\sigma$-invariant K\"ahler metric on $Y$. Let $\kappa$ be the K\"ahler 
class of $g''$. Then $\kappa$ is $\sigma$-invariant, regarded as an 
equivalence class of metrics on $Y$. By Condition \ref{n5cond} we 
know that $c_1(Y)=0$ and $h^{2,0}(Y)=0$, and that $Y\setminus S$ is 
simply-connected, where $S=\{p_1,\ldots,p_k\}$ is the singular 
set of $Y$. Thus by Proposition \ref{n3prop3}, the K\"ahler class $\kappa$ 
contains a unique metric $g_{\sst Y}$ such that $(Y,J,g_{\sst Y})$ 
is a Calabi--Yau orbifold. As $\kappa$ is $\sigma$-invariant we see that
$g_{\sst Y}$ is $\sigma$-invariant, by uniqueness of~$g_{\sst Y}$.

Proposition \ref{n3prop2} shows that there exists a holomorphic
volume form $\theta$ on $Y$. Since $\sigma$ is antiholomorphic, it
is easy to show that $\sigma^*(\theta)={\rm e}^{i\phi}\bar\theta$, for 
some $\phi\in[0,2\pi)$. Define $\theta_{\sst Y}={\rm e}^{i\phi/2}\theta$.
Then $\theta_{\sst Y}$ is a holomorphic volume form for $(Y,J,g_{\sst Y})$,
and $\sigma^*(\theta_{\sst Y})=\bar\theta_{\sst Y}$, as we want.

Let $(\Omega_{\sst Y},g_{\sst Y})$ be as in Proposition 
\ref{n3prop4}. Then $\Omega_{\sst Y}=\ha\omega_{\sst Y}\w
\omega_{\sst Y}+\Re(\theta_{\sst Y})$, where $\omega_{\sst Y}$ 
is the K\"ahler form of $g_{\sst Y}$. As 
$\sigma^*(g_{\sst Y})=g_{\sst Y}$ and $\sigma^*(J)=-J$ we have
$\sigma^*(\omega_{\sst Y})=-\omega_{\sst Y}$, and 
$\sigma^*(\Re(\theta_{\sst Y}))=\Re(\theta_{\sst Y})$ as 
$\sigma^*(\theta_{\sst Y})=\bar\theta_{\sst Y}$. Thus 
$\Omega_{\sst Y}$ and $g_{\sst Y}$ are both $\sigma$-invariant.
\end{proof}

In our next result, if $Y$ is an orbifold and $p\in Y$ an orbifold 
point modelled on $\R^n/G$, then we say that the {\it tangent space} 
$T_pY$ to $Y$ at $p$ is $\R^n/G$, in the obvious way. The proof
looks complicated, but it's really only linear algebra.

\begin{prop} For each\/ $j=1,\ldots,k$ we can identify 
the tangent space $T_{\smash{p_j}}Y$ to $Y$ at $p_j$ with\/ 
$\C^4/\an{\alpha}$ so that\/ $g_{\sst Y}$ is identified with\/ 
$\ms{\d z_1}+\cdots+\ms{\d z_4}$ at\/ $p_j$, and\/ 
$\theta_{\sst Y}$ is identified with\/ $\d z_1\w\cdots\w\d z_4$
at\/ $p_j$, and\/ $\d\sigma:T_{\smash{p_j}}Y\ra T_{\smash{p_j}}Y$ 
is identified with the map $\beta:\C^4/\an{\alpha}\ra\C^4/\an{\alpha}$ 
given by
\begin{equation}
\beta:(z_1,\dots,z_4)\an{\alpha}\longmapsto
(\bar z_2,-\bar z_1,\bar z_4,-\bar z_3)\an{\alpha}.
\label{n5beeq}
\end{equation}
\label{n5prop2}
\end{prop}

\begin{proof} Since $J,g_{\sst Y}$ and $\theta_{\sst Y}$ form a 
Calabi--Yau structure on $Y$, there certainly exists an
isomorphism $\iota:T_{\smash{p_j}}Y\ra\C^4/\an{\alpha}$ which
identifies $g_{\sst Y}$ with $\ms{\d z_1}+\cdots+\ms{\d z_4}$ 
and $\theta_{\sst Y}$ with $\d z_1\w\cdots\w\d z_4$. This $\iota$ 
is unique up to the action of $\SU(4)$ on $\C^4/\an{\alpha}$.
That is, if $B\in\SU(4)$ then $B\circ\iota:T_{\smash{p_j}}Y
\ra\C^4/\an{\alpha}$ also identifies $g_{\sst Y}$ with 
$\ms{\d z_1}+\cdots+\ms{\d z_4}$ and $\theta_{\sst Y}$ 
with~$\d z_1\w\cdots\w\d z_4$.

Now $\d\sigma:T_{\smash{p_j}}Y\ra T_{\smash{p_j}}Y$ is complex
antilinear, and so $\iota$ identifies $\d\sigma$ with the map
$\gamma:\C^4/\an{\alpha}\ra\C^4/\an{\alpha}$ given by
\begin{equation}
\!\!\!\!\!\!\!\!\!\!\!\!\!\!\!\!\!\!\!\!\!\!\!\!\!\!\!\!\!
\gamma:\left\{i^k\hbox{\small$\begin{pmatrix}z_1 \\ z_2 \\ z_3 \\ z_4
\end{pmatrix}$}:k=0,1,2,3\right\}\longmapsto\left\{i^kA
\hbox{\small$\begin{pmatrix}\bar z_1 \\ \bar z_2 \\ \bar z_3 \\ 
\bar z_4\end{pmatrix}$}:k=0,1,2,3\right\},
\!\!\!\!\!
\label{n5gaeq}
\end{equation}
for some $4\times 4$ complex matrix $A$. In fact $A$ is only defined
up to multiplication by a power of~$i$. 

As $\d\sigma$ preserves $g_{\sst Y}$ and takes $\theta_{\sst Y}$ to 
$\bar\theta_{\sst Y}$ on $T_{p_j}Y$, it follows that $\gamma$ 
preserves $\ms{\d z_1}+\cdots+\ms{\d z_4}$ and takes 
$\d z_1\w\cdots\w\d z_4$ to $\d\bar z_1\w\cdots\w\d\bar z_4$ on 
$\C^4/\an{\alpha}$. These imply that $A\bar A^t=I$ and $\det(A)=1$,
and so $A\in\SU(4)$. Also, $\gamma^2=I$ as $\sigma^2=\id$, and this
implies that $A\bar A=i^kI$ for $k=0,1,2$ or 3. And because $\sigma$
fixes only $p_1,\ldots,p_k$ in $Y$, the only fixed point of $\gamma$ 
in $\C^4/\an{\alpha}$ is~0.

So $A$ lies in $\SU(4)$ and satisfies $A\bar A=i^kI$. When we 
replace $\iota$ by $B\circ\iota$ for $B\in\SU(4)$, the matrix 
$A$ is replaced by $BAB^t$. We wish to show that we can choose 
$B\in\SU(4)$ such that the maps $\beta$ of \eq{n5beeq} and $\gamma$ 
of \eq{n5gaeq} coincide. That is, we must show that there exists 
$B\in\SU(4)$ and $l=0,1,2$ or 3 such that
\begin{equation}
i^lBAB^t=\hbox{\small
$\begin{pmatrix} 0 & 1 & 0 & 0 \\ -1 & 0 & 0 & 0 \\
0 & 0 & 0 & 1 \\ 0 & 0 & -1 & 0 \end{pmatrix}$.}
\label{ibabeq}
\end{equation}

Now $A\bar A=i^kI$ shows that $A$ and $\bar A$ commute, and so 
$A\bar A=\bar AA=\overline{A\bar A}$. Thus $i^kI$ is a real matrix,
which implies that $k=0$ or 2, and $A\bar A=\pm I$. By studying 
the eigenvectors of $A$, one can prove that there exists 
$B\in\SU(4)$ such that $BAB^t$ is one of
\begin{equation*}
I,\;\> -I,\;\>
\hbox{\small
$\begin{pmatrix} 1 & 0 & 0 & 0 \\ 0 & 1 & 0 & 0 \\
0 & 0 & -1 & 0 \\ 0 & 0 & 0 & -1 \end{pmatrix}$},\;\>
\hbox{\small$\begin{pmatrix} 0 & 1 & 0 & 0 \\ -1 & 0 & 0 & 0 \\
0 & 0 & 0 & 1 \\ 0 & 0 & -1 & 0 \end{pmatrix}$},\;\>
i\hbox{\small$\begin{pmatrix} 0 & 1 & 0 & 0 \\ -1 & 0 & 0 & 0 \\
0 & 0 & 0 & 1 \\ 0 & 0 & -1 & 0 \end{pmatrix}$.}
\end{equation*}
We exclude the first three possibilities because $\gamma$ fixes 
$(1,0,0,0)\an{\alpha}$ in $\C^4/\an{\alpha}$, contradicting the 
fact that the only fixed point of $\gamma$ in $\C^4/\an{\alpha}$ 
is 0. Putting $l=0$ in the fourth case and $l=3$ in the fifth, 
we see that \eq{ibabeq} holds. Thus $B\circ\iota$ identifies 
$T_{\smash{p_j}}Y$ with $\C^4/\an{\alpha}$ and satisfies all 
the conditions of the proposition, and the proof is complete.
\end{proof}

Now \S\ref{n41} defined a finite group $G=\an{\alpha,\beta}$ acting 
on $\R^8$, and the definitions \eq{n5aleq} and \eq{n5beeq} of 
$\alpha$ and $\beta$ above coincide with \eq{n4abeq1} in \S\ref{n41}. 
Thus the singularities of $Z=Y/\an{\sigma}$ are all modelled on 
$\R^8/G$, and we easily prove:

\begin{cor} Define $Z=Y/\an{\sigma}$. Then $Z$ is a compact, 
real\/ $8$-dimensional orbifold. The $\Spin(7)$-structure 
$(\Omega_{\sst Y},g_{\sst Y})$ on $Y$ pushes down to give a 
torsion-free $\Spin(7)$-structure $(\Omega_{\sst Z},g_{\sst Z})$ 
on $Z$. The singularities of\/ $Z$ are $k$ points $p_1,\ldots,p_k$. 
For each\/ $j=1,\ldots,k$ there is an isomorphism 
$\iota_j:\R^8/G\ra T_{\smash{p_j}}Z$ which identifies the 
$\Spin(7)$-structures $(\Omega_0,g_0)$ on $\R^8/G$ and\/ 
$(\Omega_{\sst Z},g_{\sst Z})$ on $T_{\smash{p_j}}Z$. Here $G$ 
and\/ $(\Omega_0,g_0)$ are defined in~\S\ref{n41}.
\label{n5cor}
\end{cor}

\subsection{Desingularizing $Z$ to get a compact 8-manifold $M$}
\label{n52}

So far we have constructed a $\Spin(7)$-orbifold 
$(Z,\Omega_{\sst Z},g_{\sst Z})$ with finitely many singular 
points $p_1,\ldots,p_k$, each modelled on the singularity 
$\R^8/G$ of \S\ref{n41}. But in \S\ref{n41} and \S\ref{n42} 
we wrote down two ALE $\Spin(7)$-manifolds $X_1$ and $X_2$ 
asymptotic to $\R^8/G$. We shall now resolve each singular 
point $p_j$ in $Z$ using either $X_1$ or $X_2$ to get a 
compact 8-manifold $M$. We include a parameter $t\in(0,1]$
in the construction.

\begin{dfn} For each $j$ let $\iota_j$ be as in Corollary
\ref{n5cor}, and let $\exp_{\smash{p_j}}:T_{\smash{p_j}}Z\ra Z$ 
be the {\it exponential map}, which is well-defined as $Z$ is 
complete. Then $\exp_{\smash{p_j}}\circ\,\iota_j$ maps $\R^8/G$ to $Z$.
Choose $\zeta>0$ small, and let $B_{2\zeta}(\R^8/G)$ be the open ball of 
radius $2\zeta$ about 0 in $\R^8/G$. Define $U_j\subset Z$ by 
$U_j=\exp_{\smash{p_j}}\circ\,\iota_j\bigl(B_{2\zeta}(\R^8/G)\bigr)$, 
and $\psi_j:B_{2\zeta}(\R^8/G)\ra U_j$ by $\psi_j=\exp_{\smash{p_j}}
\circ\,\iota_j$. Let $\zeta>0$ be chosen small enough that $U_j$ is 
open in $Z$ and $\psi_j:B_{2\zeta}(\R^8/G)\ra U_j$ is a diffeomorphism 
for $1\le j\le k$, and that $U_i\cap U_j=\emptyset$ when~$i\ne j$.
\label{n5zedfn}
\end{dfn}

\begin{prop} There is a smooth\/ $3$-form $\sigma_j$ on 
$B_{2\zeta}(\R^8/G)$ for $1\!\le\!j\!\le\!k$ and a constant\/ $C_1>0$, 
such that\/ $\psi_{\smash{j}}^*(\Omega_{\sst Z})-\Omega_0=\d\sigma_j$ 
and\/ $\md{\nabla^l\sigma_j}\le C_1r^{3-l}$ on $B_{2\zeta}(\R^8/G)$, 
for $l=0,1,2$. Here $\md{\,.\,}$ and\/ $\nabla$ are defined using 
the metric $g_0$ on $B_{2\zeta}(\R^8/G)$, and\/ $r:B_{2\zeta}(\R^8/G)
\ra[0,2\zeta)$ is the radius function.
\label{n5prop3}
\end{prop}

\begin{proof} The derivative of $\exp_{\smash{p_j}}$ at 0 is
the identity map on $T_{\smash{p_j}}Z$. Thus the derivative of 
$\psi_j$ at 0 is $\iota_j:\R^8/G\ra T_{\smash{p_j}}Z$, 
and so $\psi_j^*(\Omega_{\sst Z})\vert_0=\iota_j^*(\Omega_{\sst Z})=
\Omega_0\vert_0$, since $\iota_j$ identifies $\Omega_0$ and 
$\Omega_{\sst Z}$. Therefore $\psi_j^*(\Omega_{\sst Z})=\Omega_0$ at 0 
in $B_{2\zeta}(\R^8/G)$. As $\psi_j^*(\Omega_{\sst Z})-\Omega_0$ 
is a 4-form on a subset of $\R^8/G$, we can pull it back 
to $\R^8$, and regard $\psi_j^*(\Omega_{\sst Z})-\Omega_0$ as a 
4-form on the ball $B_{2\zeta}(\R^8)$ of radius $2\zeta$ in~$\R^8$. 

Then $\psi_j^*(\Omega_{\sst Z})-\Omega_0$ is a smooth 
$G$-invariant 4-form on $B_{2\zeta}(\R^8)$ which 
vanishes at 0. But $G$ contains $-1:\R^8\ra\R^8$, and any 4-form 
invariant under this map $-1$ has zero first derivative at 0. 
Hence $\psi_j^*(\Omega_{\sst Z})-\Omega_0$ vanishes to first order at
0 in $B_{2\zeta}(\R^8)$, and so by Taylor's Theorem we can show that
$\bmd{\psi_j^*(\Omega_{\sst Z})-\Omega_0}=O(r^2)$ and $\bmd{\nabla\psi_j^*
(\Omega_{\sst Z})-\Omega_0}=O(r)$ on~$B_{2\zeta}(\R^8)$.

Now $\Omega_{\sst Z}$ and $\Omega_0$ are closed, so that
$\psi_j^*(\Omega_{\sst Z})-\Omega_0$ is closed, and as 
$B_{2\zeta}(\R^8/G)$ is contractible we can write 
$\psi_j^*(\Omega_{\sst Z})-\Omega_0=\d\sigma_j$ for some smooth 3-form 
$\sigma_j$ on $B_{2\zeta}(\R^8/G)$. Since 
$\psi_j^*(\Omega_{\sst Z})-\Omega_0$ vanishes to first order at 0 we
can easily arrange that $\sigma_j$ vanishes to second order at 0,
and therefore $\md{\nabla^l\sigma}=O(r^{3-l})$ for $l=0,1,2$, using 
Taylor's Theorem as above. Thus there exists $C_1>0$ such that 
$\md{\nabla^l\sigma_j}\le C_1r^{3-l}$ on $B_{2\zeta}(\R^8/G)$, 
for $l=0,1,2$ and~$j=1,\ldots,k$.
\end{proof}

\begin{dfn} Let the ALE $\Spin(7)$-manifolds $(X_n,\Omega_n,g_n)$
and projections $\pi_n:X_n\ra\R^8/G$ be as in \S\ref{n41} and 
\S\ref{n42} for $n=1,2$. For each $t\in(0,1]$ and $n=1,2$ let 
$X_n^t=X_n$, define a $\Spin(7)$-structure $(\Omega_n^t,g_n^t)$ 
on $X_n^t$ by $\Omega_n^t=t^4\Omega_n$ and $g_n^t=t^2g_n$, and define 
$\pi_n^t:X_n^t\ra\R^8/G$ by $\pi_n^t=t\pi_n$. Then $(X_n^t,\Omega_n^t,
g_n^t)$ is an ALE $\Spin(7)$-manifold asymptotic to~$\R^8/G$.

Using the ideas of \cite{Joy6} or the explicit formula of
Calabi \cite[p.~285]{Cal} we can show that there exists 
$C_2>0$ and a smooth 3-form $\tau_n^t$ on $\R^8/G
\big\backslash B_{t\zeta}(\R^8/G)$, satisfying
\begin{equation}
\!\!\!\!\!\!\!\!\!\!\!\!\!\!\!\!\!\!\!\!\!\!\!\!\!\!\!\!\!\!\!\!\!
(\pi_n^t)_*(\Omega_n^t)=\Omega_0+\d\tau_n^t
\quad\text{and}\quad
\bmd{\nabla^l\tau_n^t}\le C_2t^8r^{-7-l}
\quad\text{for $l=0,1,2$}
\!\!\!\!\!
\label{n5tauieq}
\end{equation}
on $\R^8/G\big\backslash B_{t\zeta}(\R^8/G)$, where $\md{\,.\,}$ 
and $\nabla$ are defined using the metric~$g_0$.
\label{n5tauidef}
\end{dfn}

For $j=1,\ldots,k$, choose $n_j$ to be 1 or 2. There are $2^k$ 
ways of defining the $n_j$. We shall resolve each singular point 
$p_j$ in $Z$ using $X_{n_j}^t$ to get a 1-parameter family of 
resolutions $(M^t,\pi^t)$ of~$Z$.

\begin{dfn} For each $j=1,\ldots,k$, define open subsets $M_0^t$ 
in $Z$ and $M_j^t$ in $X_{\smash{n_j}}^t$ for $1\le j\le k$ by
\begin{equation*}
M_0^t=Z\,\Big\backslash\bigcup_{j=1}^k\psi_j
\bigr(\overline B_{t^{4/5}\zeta}(\R^8/G)\bigr)
\quad\text{and}\quad
M_j^t=(\pi_{n_j}^t)^{-1}\bigl(B_{2t^{4/5}\zeta}(\R^8/G)\bigr).
\end{equation*}
That is, $M_0^t$ is the complement in $Z$ of the closed balls of 
radius $t^{4/5}\zeta$ about $p_j$ for $1\le j\le k$, and $M_j^t$ is 
the inverse image of $B_{2t^{4/5}\zeta}(\R^8/G)$ in~$X_{\smash{n_j}}^t$.

Define an equivalence relation `$\sim$' on the disjoint union
$\coprod_{j=0}^kM_j^t$ by $x\sim y$ if either (a) $x=y$,
\begin{itemize}
\item[(b)] $x\in M_j^t$ and $y\in U_j\cap M_0^t$ and 
$\psi_j\circ\pi_{n_j}^t(x)=y$, for some $j=1,\ldots,k$, or
\item[(c)] $y\in M_j^t$ and $x\in U_j\cap M_0^t$ and 
$\psi_j\circ\pi_{n_j}^t(y)=x$, for some $j=1,\ldots,k$.
\end{itemize}
Define the {\it resolution} $M^t$ of $Z$ to be
$\coprod_{j=0}^kM_j^t/\sim$. It is easy to see that $M^t$
is a compact 8-manifold. Define a projection $\pi^t:M^t\ra Z$
by $\pi^t\bigl([x]\bigr)=x$ when $x\in M_0^t$, and 
$\pi^t\bigl([x]\bigr)=\psi_j\circ\pi_{n_j}^t(x)=x$ when
$x\in M_j^t$ for some $j=1,\ldots,k$, where $[x]$ is the
equivalence class of $x$ under $\sim$. Then $\pi^t$ is
well-defined, continuous and surjective, and 
$\pi^t:M^t\big\backslash\bigcup_{j=1}^k(\pi^t)^{-1}(p_j)\ra
Z\,\big\backslash\{p_1,\ldots,p_k\}$ is a diffeomorphism.
\label{n5resdef}
\end{dfn}

Since the resolutions $(M^t,\pi^t)$ of $Z$ form a smooth connected
family, they are all diffeomorphic to the same compact 8-manifold $M$. 
We can regard $M_j^t$ as an open subset of $M^t$ for $j=0,\ldots,k$,
and then the $M_j^t$ form an {\it open cover} of $M^t$. If 
$1\le i,j\le k$ and $i\ne j$ then $M_i^t\cap M_j^t=\emptyset$.
The overlap $M_0^t\cap M_j^t$ is naturally isomorphic to an 
{\it annulus} in $\R^8/G$, with inner radius 
$t^{4/5}\zeta$ and outer radius $2t^{4/5}\zeta$. The reason for 
including the factors $t^{4/5}$ will be explained shortly.

We now calculate the fundamental group of $M^t$.

\begin{prop} If\/ $n_j=1$ for $j=1,\ldots,k$ then 
$\pi_1(M^t)\cong\Z_2$. Otherwise, $M^t$ is simply-connected.
\label{n5prop4}
\end{prop}

\begin{proof} Since $Y\setminus\{p_1,\ldots,p_k\}$ is 
simply-connected by Condition \ref{n5cond} and $\sigma$ acts 
freely on $Y\setminus\{p_1,\ldots,p_k\}$, we see that the 
fundamental group of $Z\setminus\{p_1,\ldots,p_k\}$ is $\Z_2$. The 
natural inclusion of $Z\setminus\{p_1,\ldots,p_k\}$ in $M^t$ induces 
a homomorphism from $\pi_1\bigl(Z\setminus\{p_1,\ldots,p_k\}\bigr)$ 
to $\pi_1(M^t)$, which is easily shown to be surjective. Also, as 
$X_{n_j}^t$ is $X_1$ or $X_2$ we have~$\pi_1(X_{n_j}^t)\cong\Z_2$. 

Therefore, $\pi_1(M_t)$ is $\Z_2$ if the generator of 
$\pi_1\bigl(Z\setminus\{p_1,\ldots,p_k\}\bigr)$ projects to the 
nonzero element of $\pi_1(X_{n_j}^t)$ for all $1\le j\le k$, and 
$\pi_1(M^t)$ is trivial otherwise. But calculation shows that 
the generator of $\pi_1\bigl(Z\setminus\{p_1,\ldots,p_k\}\bigr)$ 
is nonzero in $\pi_1(X_{n_j}^t)$ if and only if~$n_j=1$.
\end{proof}

This shows that of the $2^k$ possible ways of choosing 
the $n_j$, one possibility gives $\pi_1(M^t)=\Z_2$, and the 
remaining $2^k-1$ possibilities all give simply-connected~$M^t$. 

\subsection{A $\Spin(7)$-structure $(\Omega^t,g^t)$ on $M^t$ with 
small torsion}
\label{n53}

Each open subset $M_j^t$ in $M^t$ carries a torsion-free 
$\Spin(7)$-structure, $(\Omega_{\sst Z},g_{\sst Z})$ for $j=0$
and $(\Omega_{\smash{n_j}}^t,g_{\smash{n_j}}^t)$ for $1\le j\le k$. 
We shall join these $\Spin(7)$-structures together with a partition 
of unity to get a $\Spin(7)$-structure $(\Omega^t,g^t)$ on $M^t$ and 
estimate its torsion.

\begin{dfn} Let $\eta:[0,\infty)\ra[0,1]$ be a smooth function with 
$\eta(x)=0$ for $x\le\zeta$ and $\eta(x)=1$ for $x\ge 2\zeta$. Define 
a 4-form $\xi^t$ on $M^t$ by $\xi^t=\Omega_{\sst Z}$ in $M_0^t\big
\backslash\bigcup_{j=1}^kM_j^t$, and $\xi^t=\Omega_{\smash{n_j}}^t$ 
in $M_j^t\setminus M_0^t$ for $1\le j\le k$, and
\begin{equation}
\!\!\!\!\!\!\!\!\!\!\!\!\!\!\!\!\!\!\!\!\!\!\!\!\!\!\!\!
\xi^t=\Omega_0+\d\bigl(\eta(t^{-4/5}r)\sigma_j\bigr)
+\d\bigl((1-\eta(t^{-4/5}r))\tau_{n_j}^t\bigr)
\quad\text{in $M_0^t\cap M_j^t$} 
\!\!\!\!
\label{n5xiteq}
\end{equation}
for $1\le j\le k$, where  we identify $M_0^t\cap M_j^t$ with an 
annulus in $\R^8/G$ in the natural way. Since $\Omega_{\sst Z}=\Omega_0+
\d\sigma_j$ and $\Omega_{\smash{n_j}}^t=\Omega_0+\d\tau_{\smash{n_j}}^t$
in $M_0^t\cap M_j^t$ it follows that $\xi^t$ is smooth, and as 
$\Omega_{\sst Z}$, $\Omega_{\smash{n_j}}^t$ and $\Omega_0$ are closed, 
$\xi^t$ is closed.
\label{n5xitdef}
\end{dfn}

\begin{lem} There exists $C_3>0$ such that for each\/ $j=1,\ldots,k$
and\/ $t\in(0,1]$, this $4$-form $\xi^t$ satisfies
\begin{equation}
\bmd{\xi^t-\Omega_0}\le C_3t^{8/5}
\quad\text{and}\quad
\bmd{\nabla(\xi^t-\Omega_0)}\le C_3t^{4/5}
\label{n5xitesteq}
\end{equation}
in $M_0^t\cap M_j^t$, where $\md{\,.\,}$ and\/ $\nabla$ are
defined using the metric~$g_0$.
\label{n5lem}
\end{lem}

\begin{proof} Expanding \eq{n5xiteq} we find that
\begin{equation*}
\xi^t-\Omega_0=\eta(t^{-4/5}r)\d\sigma_j+
(1-\eta(t^{-4/5}r))\d\tau_{n_j}^t
+t^{-4/5}\eta'(t^{-4/5}r)\d r\w(\sigma_j-\tau_{n_j}^t)
\end{equation*}
in $M_0^t\cap M_j^t$. Since $t^{4/5}\zeta\le r\le 2t^{4/5}\zeta$,
Proposition \ref{n5prop3} and \eq{n5tauieq} show that
\begin{alignat*}{3}
\bmd{\sigma_j}&\le 8C_1\zeta^3t^{12/5},\quad &
\bmd{\d\sigma_j}&\le 4C_1\zeta^2t^{8/5},\quad &
\bmd{\nabla\d\sigma_j}&\le 2C_1\zeta t^{4/5},\\
\bmd{\tau_{n_j}^t}&\le C_2\zeta^{-7}t^{12/5},\quad &
\bmd{\d\tau_{n_j}^t}&\le C_2\zeta^{-8}t^{8/5},\quad & \text{and}\quad
\bmd{\nabla\d\tau_{n_j}^t}&\le C_2\zeta^{-9}t^{4/5}.
\end{alignat*}
Combining these with the previous equation and using the facts
that $\md{\d r}=1$ and $\eta'$ is bounded independently of
$t$, we soon prove~\eq{n5xitesteq}.
\end{proof}

We can now explain why we chose the power $t^{4/5}$ in Definition
\ref{n5resdef}. Suppose we had defined $M^t$ and $\xi^t$ using
$t^\alpha$ in place of $t^{4/5}$, for some $\alpha\in[0,1]$. 
Then in the calculation above the $\sigma_j$ and 
$\tau_{\smash{n_j}}^t$ terms would contribute $O(t^{2\alpha})$ 
and $O(t^{8-8\alpha})$ to $\xi^t-\Omega_0$ respectively, and so 
$\xi^t-\Omega_0$ would be $O(t^{2\alpha})+O(t^{8-8\alpha})$. 
This is smallest when $2\alpha=8-8\alpha$, that is, when 
$\alpha=4/5$. So the power $t^{4/5}$ minimizes the size 
of~$\xi^t-\Omega_0$.

Now we can define the $\Spin(7)$-structures $(\Omega^t,g^t)$ on~$M^t$.

\begin{dfn} Let\/ $\rho$ be as in Proposition \ref{n2prop},
and choose $\epsilon\in(0,1]$ such that\/ $C_3\epsilon^{8/5}\le\rho$.
Suppose $t\in(0,\epsilon]$. Then $\bmd{\xi^t-\Omega_0}\le C_3t^{8/5}\le
\rho$ in $M_0^t\cap M_j^t$ for $1\le j\le k$ by \eq{n5xitesteq},
and so $\xi^t$ lies in ${\mathcal T}M^t$ on $M_0^t\cap M_j^t$ by 
part (i) of Proposition \ref{n2prop}. But $\xi^t$ is $\Omega_{\sst Z}$ 
or $\Omega_{\smash{n_j}}^t$ outside the overlaps $M_0^t\cap M_j^t$, and 
thus $\xi^t\in C^\infty({\mathcal T}M^t)$. For each $t\in(0,\epsilon]$ 
define $\Omega^t=\Theta(\xi^t)$, where $\Theta$ is given in Proposition
\ref{n2prop}. Then $\Omega^t\in C^\infty({\mathcal A}M^t)$, and so 
$\Omega^t$ extends to a $\Spin(7)$-structure $(\Omega^t,g^t)$ on $M^t$. 
Define a 4-form $\phi^t$ on $M^t$ by $\phi^t=\xi^t-\Omega^t$. 
Then $\d\Omega^t+\d\phi^t=0$, as $\d\xi^t=0$ on~$M^t$. 
\label{n5Omtdef}
\end{dfn}

Here $\xi^t$ is a 4-form which does not lie in ${\mathcal A}M^t$, 
but is close to ${\mathcal A}M^t$ for small $t$, and $\Omega^t$ is 
the section of ${\mathcal A}M^t$ closest to $\xi^t$. What is really 
happening is that the $\Spin(7)$-structure $(\Omega^t,g^t)$ is equal to 
$(\Omega_{\smash{n_j}}^t,g_{\smash{n_j}}^t)$ in $M_j^t\setminus M_0^t$
and to $(\Omega_{\sst Z},g_{\sst Z})$ outside $M_j^t$ for $j=1,\ldots,k$,
and $(\Omega^t,g^t)$ interpolates smoothly between these two possibilities
on the annulus~$M_j^t\cap M_0^t$.

\subsection{Existence of torsion-free $\Spin(7)$-structures on $M$}
\label{n54}

Next we shall show that $(\Omega^t,g^t)$ can be deformed to a
torsion-free $\Spin(7)$-structure on $M$ when $t$ is small.

\begin{thm} In the situation above, there exist constants 
$\lambda,\mu,\nu>0$ such that for all\/ $t\in(0,\epsilon]$ we have
\begin{itemize}
\item[{\rm(i)}] $\lnm{\phi^t}2\le\lambda t^{24/5}$ 
and\/~$\lnm{\d\phi^t}{10}\le\lambda t^{36/25}$;
\item[{\rm(ii)}] the injectivity radius $\delta(g^t)$ satisfies
$\delta(g^t)\ge\mu t$; and
\item[{\rm(ii)}] the Riemann curvature $R(g^t)$ 
satisfies~$\bcnm{R(g^t)}0\le\nu t^{-2}$.
\end{itemize}
Here all norms are calculated using the metric $g^t$ on~$M^t$.
\label{n5thm1}
\end{thm}

\begin{proof} Outside the overlaps $M_0^t\cap M_j^t$ for
$1\le j\le k$ we either have $\xi^t=\Omega^t=\Omega_{\sst Z}$ or 
$\xi^t=\Omega^t=\Omega_{\smash{n_j}}^t$. In both cases 
$\phi^t=\xi^t-\Omega^t=0$, and so $\phi^t$ is zero outside the 
$M_0^t\cap M_j^t$. In $M_0^t\cap M_j^t$ we apply part (ii) of 
Proposition \ref{n2prop} with $\Omega=\Omega_0$ and $\xi=\xi^t$, to get
\begin{equation*}
\bmd{\phi^t}_{g^t}\le \md{\xi^t-\Omega_0}_{g_0}
\quad\text{and}\quad
\bmd{\nabla^{g^t}\phi^t}_{g^t}\le C
\bmd{\nabla^{g_0}(\xi^t-\Omega_0)}_{g_0}.
\end{equation*}
Combining this with \eq{n5xitesteq} gives
\begin{equation*}
\bmd{\phi^t}_{g^t}\le C_3t^{8/5}
\quad\text{and}\quad
\bmd{\d\phi^t}_{g^t}\le\bmd{\nabla^{g^t}\phi^t}_{g^t}\le CC_3t^{4/5}.
\end{equation*}

Now each $M_0^t\cap M_j^t$ is an annulus in $\R^8/G$ with inner 
radius $t^{4/5}\zeta$ and outer radius $2t^{4/5}\zeta$, and the metric 
$g^t$ on $M_0^t\cap M_j^t$ is close to the flat metric $g_0$ on 
$\R^8/G$. Therefore we can find $C_4>0$ independent of $t$ such 
that $\sum_{j=1}^k\vol\bigl(M_0^t\cap M_j^t)\le C_4t^{32/5}$. Hence
\begin{equation*}
\int_{M^t}\ms{\phi^t}\d V\le (C_3t^{8/5})^2C_4t^{32/5}
\quad\text{and}\quad
\int_{M^t}\md{\d\phi^t}^{10}\d V\le (CC_3t^{4/5})^{10}C_4t^{32/5}.
\end{equation*}
Taking roots gives part (i) of the theorem, with
$\lambda=C_3\max(C_4^{1/2},CC_4^{1/10})$. 

Parts (ii) and (iii) are elementary. The metric $g_{\smash{n_j}}^t$ 
is made by scaling $g_{\smash{n_j}}$ by a factor $t$. Thus 
$\delta(g_{\smash{n_j}}^t)=t\delta(g_{\smash{n_j}})$ and 
$\cnm{R(g_{\smash{n_j}}^t)}0=t^{-2}\cnm{R(g_{\smash{n_j}})}0$. 
We make $g^t$ by gluing together the $g_{\smash{n_j}}^t$ on the 
patches $M_j^t$ for $j=1,\ldots,k$ and $g_{\sst Z}$ on $M_0^t$.
It is clear that for small $t$, the dominant contributions to 
$\delta(g^t)$ and $\cnm{R(g^t)}0$ come from $\delta(g_{\smash{n_j}}^t)$ 
and $\cnm{R(g_{\smash{n_j}}^t)}0$ for some $j$, and these are 
proportional to $t$ and $t^{-2}$. This proves (ii) and (iii) 
for some $\mu,\nu>0$, and the theorem is complete.
\end{proof}

Finally we can prove our main result.

\begin{thm} Suppose Condition \ref{n5cond} holds, and let\/ $M$ be 
the compact\/ $8$-manifold defined in Definition \ref{n5resdef}. Then 
there exist torsion-free $\Spin(7)$-structures $(\tilde\Omega,\tilde g)$ 
on $M$. If\/ $\pi_1(M)=\{1\}$ then $\Hol(\tilde g)=\Spin(7)$, and if\/ 
$\pi_1(M)=\Z_2$ then~$\Hol(\tilde g)=\Z_2\ltimes\SU(4)$.
\label{n5thm2}
\end{thm}

\begin{proof} Let $\lambda,\mu,\nu$ be as in Theorem \ref{n5thm1}. 
Then Theorem \ref{n2thm3} gives a constant $\kappa>0$. Choose 
$t>0$ with $t\le\epsilon\le 1$ and $t\le\kappa$. Let $(\Omega,g)$ be the 
$\Spin(7)$-structure $(\Omega^t,g^t)$ on $M=M^t$, and $\phi$ the 
4-form $\phi^t$. Then $\d\Omega+\d\phi=0$ by Definition \ref{n5Omtdef}, 
and parts (i)--(iii) of Theorem \ref{n5thm1} imply (i)--(iii) 
of Theorem \ref{n2thm3}, as~$t\le 1$.

Therefore all the hypotheses of Theorem \ref{n2thm3} 
hold, and the theorem shows that there exists a torsion-free 
$\Spin(7)$-structure $(\tilde\Omega,\tilde g)$ on $M$. It remains 
to identify the holonomy group $\Hol(\tilde g)$ of $\tilde g$. We
can regard the $\Spin(7)$-orbifold $(Z,\Omega_{\sst Z},g_{\sst Z})$
as the limit as $t\ra 0$ of the $\Spin(7)$-manifolds
$(M,\tilde\Omega,\tilde g)$. Because of this, it is not difficult 
to show that~$\Hol(g_{\sst Z})\subseteq\Hol(\tilde g)$.

Now $\Hol(g_{\sst Z})=\Z_2\ltimes\SU(4)$, and thus
$\Z_2\ltimes\SU(4)\subseteq\Hol(\tilde g)\subseteq\Spin(7)$. If
$\pi_1(M)=\{1\}$ then $\Hol(\tilde g)$ is connected.
But the only connected Lie subgroup of $\Spin(7)$ containing
$\Z_2\ltimes\SU(4)$ is $\Spin(7)$, so $\Hol(\tilde g)=\Spin(7)$. 
If $\pi_1(M)=\Z_2$ then $\Hol(\tilde g)\ne\Spin(7)$ 
by Theorem \ref{n2thm1}. This forces $\Hol^0(\tilde g)=\SU(4)$, and 
it is then easy to see that~$\Hol(\tilde g)=\Z_2\ltimes\SU(4)$.
\end{proof}

Since by Proposition \ref{n5prop4} we can always choose the $n_j$
so that $M$ is simply-connected, we can always arrange for $\tilde g$ 
to have holonomy $\Spin(7)$. When $\pi_1(M)=\Z_2$, the complex 
orbifold $Y$ has a crepant resolution $\tilde Y$, which admits 
K\"ahler metrics $\tilde g$ with holonomy $\SU(4)$, making it into 
a Calabi--Yau manifold. The action of $\sigma$ on $Y$ lifts to a 
{\it free} action of $\sigma$ on $\tilde Y$, and so 
$M=\tilde Y/\an{\sigma}$ is a compact 8-manifold. If we choose 
$\tilde g$ to be $\sigma$-invariant then it pushes down to $M$, 
and has holonomy~$\Z_2\ltimes\SU(4)$.

\section{How to apply the construction}
\label{n6}

We now explain ways of finding orbifolds $Y$ and involutions 
$\sigma:Y\ra Y$ satisfying Condition \ref{n5cond}, and how to 
calculate the Betti numbers of the resulting 8-manifolds $M$ 
with holonomy~$\Spin(7)$.

\subsection{Finding suitable Calabi--Yau 4-orbifolds $Y$}
\label{n61}

To apply the construction of \S\ref{n6} we need a source 
of compact K\"ahler 4-orbifolds $Y$ with $c_1(Y)=0$ and isolated
singularities modelled on $\C^4/\Z_4$. Fortunately, physicists
and algebraic geometers have been studying Calabi--Yau manifolds
for many years, mainly in complex dimension 3. Several powerful
methods have been developed for constructing Calabi--Yau
manifolds, and we will adapt some of these to our problem.

The main idea we shall use is borrowed from Candelas, Lynker 
and Schrimmrigk \cite{CLS}, who constructed a large number
of Calabi--Yau 3-folds as crepant resolutions of hypersurfaces
in weighted projective spaces $\CP^4_{a_0,\ldots,a_4}$. We shall 
explain their methods, beginning with {\it weighted projective 
spaces}, which are an important class of complex orbifolds.

\begin{dfn} Let $m\ge 1$ be an integer, and $a_0,a_1,\ldots,a_m$
positive integers with highest common factor 1. Let $\C^{m+1}$ have 
complex coordinates on $(z_0,\ldots,z_m)$, and define an action of 
the complex Lie group $\C^*$ on $\C^{m+1}$ by
\begin{equation}
(z_0,\ldots,z_m)\,{\buildrel u\over\longmapsto}
(u^{a_0}z_0,\ldots,u^{a_m}z_m),\qquad\text{for $u\in\C^*$.}
\label{wpsacteq}
\end{equation}
Define the {\it weighted projective space} $\CP^m_{a_0,\ldots,a_m}$
to be $\bigl(\C^{m+1}\setminus\{0\}\bigr)/\C^*$, where $\C^*$ acts on 
$\C^{m+1}\setminus\{0\}$ with the action \eq{wpsacteq}. Then 
$\CP^m_{a_0,\ldots,a_m}$ is compact and Hausdorff, and has the 
structure of a {\it complex orbifold}.

Let $[z_0,\ldots,z_m]$ be a point in $\CP^m_{a_0,\ldots,a_m}$, 
and let $k$ be the highest common factor of the set of those 
$a_j$ for which $z_j\ne 0$. If $k=1$ then $[z_0,\ldots,z_m]$ is a 
nonsingular point of $\CP^m_{a_0,\ldots,a_m}$, and if $k>1$ then 
$[z_0,\ldots,z_m]$ is an orbifold point with orbifold group~$\Z_k$.
\label{wpsdef}
\end{dfn}

We call a polynomial $f(z_0,\ldots,z_m)$ {\it weighted 
homogeneous of degree} $d$ if
\begin{equation*}
f(u^{a_0}z_0,\ldots,u^{a_m}z_m)=u^df(z_0,\ldots,z_m)
\quad\text{for all $u,z_0,\ldots,z_m\in\C$.}
\end{equation*}
Let $f$ be such a polynomial, and define a hypersurface $Y$ in 
$\CP^m_{a_0,\ldots,a_m}$ by
\begin{equation*}
Y=\bigl\{[z_0,\ldots,z_m]\in\CP^m_{a_0,\ldots,a_m}:
f(z_0,\ldots,z_m)=0\bigr\}.
\end{equation*}
Then we call $Y$ a {\it hypersurface of degree} $d$ 
in~$\CP^m_{a_0,\ldots,a_m}$. 

We say that $f$ is {\it transverse} if $f(z_0,\ldots,z_m)=0$ 
and $\d f(z_0,\ldots,z_m)=0$ have no common solutions in 
$\C^{m+1}\setminus\{0\}$. If $f$ is transverse then the only 
singular points of $Y$ are also singular points of 
$\CP^m_{a_0,\ldots,a_m}$, and $Y$ is an {\it orbifold}, 
all of whose orbifold groups are cyclic. Note that for given 
weights $a_0,\ldots,a_m$ and degree $d$, there may not exist 
any transverse polynomials~$f$. 

So let $Y$ be a hypersurface of degree $d$ in 
$\CP^m_{a_0,\ldots,a_m}$, defined by a transverse polynomial.
Using the {\it adjunction formula}, we find that $c_1(Y)=0$
if and only if $d=a_0+\cdots+a_m$. In this case it is easy to
show that $Y$ is a {\it Calabi--Yau orbifold}. Candelas et al.\
\cite{CLS} considered the case $m=4$, and used a computer to
search for Calabi--Yau 3-orbifolds of this kind, finding some
6000 examples. They then resolved the singularities of each 
to get a Calabi--Yau 3-manifold.

As we are interested in Calabi--Yau 4-orbifolds, we shall
consider hypersurfaces $Y$ in $\CP^5_{a_0,\ldots,a_5}$. 
Here is a simple class of such~$Y$. 

\begin{ex} Let $a_0,\ldots,a_5$ be positive integers with
highest common factor $\hcf(a_0,\ldots,a_5)=1$, and let 
$d=a_0+\cdots+a_5$. Usually we order the $a_j$ with 
$a_0\le a_1\le\cdots\le a_5$. Suppose that $a_j$ divides $d$ for 
$j=0,\ldots,5$, and define $k_j=d/a_j$. Define a hypersurface $Y$ 
in $\CP^5_{a_0,\ldots,a_5}$ by
\begin{equation*}
Y=\bigl\{[z_0,\ldots,z_5]\in\CP^5_{a_0,\ldots,a_5}:
z_0^{k_0}+\cdots+z_5^{k_5}=0\bigr\}.
\end{equation*}
Since $a_jk_j=d$ we see that $z_0^{k_0}+\cdots+z_5^{k_5}$ is a 
weighted homogeneous polynomial of degree $d$, and it is also 
transverse. 

Therefore $Y$ is a complex orbifold, with singularities only at the 
intersection of $Y$ with the singular set of $\CP^5_{a_0,\ldots,a_5}$. 
Since the degree $d$ of $Y$ satisfies $d=a_0+\cdots+a_5$, we have 
$c_1(Y)=0$. Also $Y$ admits K\"ahler metrics, as $\CP^5_{a_0,\ldots,a_5}$ 
is K\"ahler. So $Y$ is a compact complex orbifold with $c_1(Y)=0$, 
admitting K\"ahler metrics.
\label{n6ex1}
\end{ex}

Now to apply the construction of \S\ref{n5} the singular points 
of $Y$ must satisfy Condition \ref{n5cond}. This is a strong 
restriction on $a_0,\ldots,a_5$, which admits only a few solutions. 
However, we can get many other suitable orbifolds $Y$ by generalizing 
our construction a bit. Here are four ways to do this.

\begin{itemize}
\item {\bf Defining $Y$ by a different polynomial.} We could 
define $Y$ using some more general transverse weighted homogenous 
polynomial of degree $d$ in $z_0,\ldots,z_5$, instead of 
$z_0^{\smash{k_0}}+\cdots+z_5^{\smash{k_5}}$. The requirement 
that $a_j$ divides $d$ for $j=0,\ldots,5$ is then replaced by
some other condition on the $a_j$ and~$d$.

\item {\bf Dividing by a finite group.} Let $W$ be a Calabi--Yau 
hypersurface in $\CP^5_{a_0,\ldots,a_5}$, and $G$ a finite group 
acting on $W$ preserving its Calabi--Yau structure. Then $Y=W/G$ 
is a Calabi--Yau orbifold.

\item {\bf Partial crepant resolutions.} Let $W$ be a Calabi--Yau
hypersurface in $\CP^5_{a_0,\ldots,a_5}$ which has some 
singularities of the kind we want, together with other 
singularities that we don't want. We let $Y$ be a partial 
crepant resolution of $W$, which resolves the singularities 
that we don't want, leaving those that we do.

\item {\bf Complete intersections in $\CP^m_{a_0,\ldots,a_5}$.}
Rather than a hypersurface in $\CP^5_{a_0,\ldots,a_5}$, we take 
$Y$ to be a {\it complete intersection} of $m-4$ hypersurfaces
in $\CP^m_{a_0,\ldots,a_m}$, for some~$m>5$.
\end{itemize}

We can also use combinations of these four techniques --- for 
instance, we can take $Y$ to be a partial crepant resolution of 
$W/G$, where $W$ is a hypersurface in $\CP^5_{a_0,\ldots,a_5}$, 
and $G$ a finite group acting on~$W$.

\subsection{Antiholomorphic maps $\sigma:Y\ra Y$}
\label{n62}

Suppose we have chosen an orbifold $Y$ as above, with
isolated singular points $p_1,\ldots,p_k$. The next ingredient 
in our construction is an antiholomorphic involution $\sigma:Y\ra Y$, 
which should fix only $p_1,\ldots,p_k$. For example, suppose 
$Y$ is a hypersurface in $\CP^5_{a_0,\ldots,a_5}$. Then to 
find $\sigma$ we would look for an antiholomorphic involution 
$\sigma:\CP^5_{a_0,\ldots,a_5}\ra\CP^5_{a_0,\ldots,a_5}$ with 
$\sigma(Y)=Y$, and restrict $\sigma$ to~$Y$. 

The most obvious such $\sigma$ maps $[z_0,\ldots,z_5]\mapsto
[\bar z_0,\ldots,\bar z_5]$. But this will not do, as its 
fixed points are not isolated in $Y$. To get isolated 
fixed points we need to try something more subtle.
Here is an example of the kind of thing we mean. 

\begin{ex} In the situation of Example \ref{n6ex1}, suppose that
$a_0,\ldots,a_3$ are odd and $a_4,a_5$ even with $a_0=a_1$, $a_2=a_3$ and 
$a_4=a_5$. Define $\sigma:\CP^5_{a_0,\ldots,a_5}\ra\CP^5_{a_0,\ldots,a_5}$ 
by
\begin{equation*}
\sigma:[z_0,\ldots,z_5]\mapsto[\bar z_1,-\bar z_0,\bar z_3,
-\bar z_2,\bar z_5,\bar z_4].
\end{equation*}
As $\sigma$ swaps the pairs $z_0,z_1$ and $z_2,z_3$ and $z_4,z_5$,
we need $a_0=a_1$, $a_2=a_3$ and $a_4=a_5$ for $\sigma$ to be 
well-defined. Clearly $\sigma$ is antiholomorphic, and~$\sigma(Y)=Y$. 

Now $\sigma^2$ acts by
\begin{equation*}
\sigma^2:[z_0,\ldots,z_5]\mapsto[-z_0,-z_1,-z_2,-z_3,z_4,z_5]
\end{equation*}
But putting $u=-1$ in \eq{wpsacteq} gives $[-z_0,-z_1,-z_2,-z_3,z_4,z_5]
=[z_0,\ldots,z_5]$, as $a_0,\ldots,a_3$ are odd and $a_4,a_5$ even. 
Thus $\sigma^2=1$, and $\sigma:Y\ra Y$ is an antiholomorphic involution.

It is not difficult to show that the fixed points of $\sigma$ in
$\CP^5_{a_0,\ldots,a_5}$ are
\begin{equation*}
\bigl\{[0,0,0,0,1,{\rm e}^{i\theta}]\in\CP^5_{a_0,\ldots,a_5}:
\theta\in[0,2\pi)\bigr\}.
\end{equation*}
Now $[0,0,0,0,1,{\rm e}^{i\theta}]$ lies in $Y$ if $1+e^{k_5i\theta}=0$.
It follows that the fixed points of $\sigma$ in $Y$ are the $k_5$
isolated points
\begin{equation*}
\bigl\{[0,0,0,0,1,{\rm e}^{\pi i(2j-1)/k_5}]:j=1,2,\ldots,k_5\bigr\}.
\end{equation*}
\label{n6ex2}
\end{ex}

Observe the trick we have used here: if $a_j=a_{j+1}$ then we 
can choose $\sigma$ to act on the coordinates $z_j,z_{j+1}$ by
$(z_j,z_{j+1})\mapsto(\bar z_{j+1},-\bar z_j)$. All the fixed
points of $\sigma$ will then satisfy $z_j=z_{j+1}=0$. By doing this
with two pairs of coordinates, say $z_0,z_1$ and $z_2,z_3$, the
fixed points of $\sigma$ satisfy $z_0=z_1=z_2=z_3=0$. Thus they will
be of complex codimension 4 in $Y$, and will be {\it isolated},
as we want. 

This trick can also be adapted to more general situations, in
which $Y$ is a quotient by a finite group, or a partial crepant 
resolution, and so on. Note that as $\sigma^2$ maps $(z_j,z_{j+1})
\mapsto(-z_j,-z_{j+1})$, care must be taken to ensure that~$\sigma^2=1$.

\subsection{Calculating the Euler characteristic of\/ $Y$}
\label{n63}

To determine the Betti numbers of the 8-manifold $M$ that we 
construct, we will need to know the {\it Euler characteristic}
of $Y$. Now there are two different notions of the Euler 
characteristic of an orbifold, defined by Satake 
\cite[\S 3.3]{Sat}. The version we are interested in is the 
{\it ordinary Euler characteristic} $\chi(Y)$, which is an
integer and satisfies $\chi(Y)=\sum_{j=0}^{2n}(-1)^jb^j(Y)$. 
There is also the {\it orbifold Euler characteristic} $\chi_V(Y)$, 
which is a rational number that crops up naturally in problems
involving characteristic classes.

In the next example we explain an elementary and fairly 
crude method for finding $\chi(Y)$ in the case that $Y$ is a 
hypersurface in $\CP^m_{a_0,\ldots,a_m}$, of the kind considered 
in Example \ref{n6ex1}. It is also possible to calculate 
$\chi_V(Y)$ using Chern classes and get $\chi(Y)$ by adding on
contributions from the singular set (see for instance Candelas 
et al.\ \cite[\S 3]{CLS}), but we will not discuss this.

\begin{ex} Let $a_0,\ldots,a_m$, $k_0,\ldots,k_m$ and $d$ be 
positive integers with $a_jk_j=d$ for $j=0,\ldots,m$. For each 
$j=0,\ldots,m$, define $Y_j\subset\CP^j_{a_0,\ldots,a_j}$ by
\begin{equation*}
Y_j=\bigl\{[z_0,\ldots,z_j]\in\CP^j_{a_0,\ldots,a_j}:
z_0^{k_0}+\cdots+z_j^{k_j}=0\bigr\},
\end{equation*}
and define $\pi_j:Y_j\ra\CP^{j-1}_{a_0,\ldots,a_{j-1}}$ 
by~$\pi:[z_0,\ldots,z_j]\mapsto[z_0,\ldots,z_{j-1}]$. 

Suppose for simplicity that $a_i$ divides $a_j$ for $0\le i<j\le m$.
Then for each $j$, $\pi_j$ is a $k_j$-fold branched cover of 
$\CP^{j-1}_{a_0,\ldots,a_{j-1}}$, branched over $Y_{j-1}$. That 
is, if $p\in\CP^{j-1}_{a_0,\ldots,a_{j-1}}$ then $\pi^{-1}(p)$ is 
one point when $p\in Y_{j-1}$ and $k_j$ points when $p\notin Y_{j-1}$. 
It follows that
\begin{equation}
\begin{split}
\chi(Y_j)&=k_j\cdot\chi(\CP^{j-1}_{a_0,\ldots,a_{j-1}})+
(1-k_j)\chi(Y_{j-1})\\
&=k_jj+(1-k_j)\chi(Y_{j-1}),
\end{split}
\label{chiyjeq}
\end{equation}
since $\chi(\CP^{j-1}_{a_0,\ldots,a_{j-1}})=j$. This equation gives
$\chi(Y_j)$ in terms of $\chi(Y_{j-1})$. Hence by induction we can
write $\chi(Y_m)$ in terms of $\chi(Y_0)$. But $Y_0=\emptyset$ so
that $\chi(Y_0)=0$, and thus we determine~$\chi(Y_m)$.

If $a_i$ does not divide $a_j$ for some $0\le i<j\le m$, then $\pi_j$ is 
also branched over other parts of $\CP^{j-1}_{a_0,\ldots,a_{j-1}}$. Let 
$p=[z_0,\ldots,z_{j-1}]$ be in $\CP^{j-1}_{a_0,\ldots,a_{j-1}}\setminus 
Y_{j-1}$, and let $I$ be the set of $i$ in $\{0,\ldots,j-1\}$ for 
which $z_i\ne 0$. Define $l=\hcf(a_i:i\in I)$ and $m=\hcf(l,a_j)$. 
Then it turns out that $\pi^{-1}(p)$ is $k_jm/l$ points in $Y_j$.
Clearly $k_jm/l=k_j$ if $l=m$, that is, if $l$ divides~$a_j$.

Thus $\pi_j$ is also branched over subsets of 
$\CP^{j-1}_{a_0,\ldots,a_{j-1}}\setminus Y_{j-1}$ corresponding to 
subsets $I\subseteq\{0,\ldots,j-1\}$ for which $l=\hcf(a_i:i\in I)$
does not divide $a_j$. To calculate $\chi(Y_j)$ in this case we 
must modify \eq{chiyjeq} by adding in contributions from each 
such $I$. We will explain this when we meet it in examples later.
\label{n6ex3}
\end{ex}

\subsection{How to find topological invariants of $Y$, $Z$ and $M$}
\label{n64}

To calculate the cohomology and fundamental group of our complex
orbifolds $Y$ we will need the following result, a form of the 
{\it Lefschetz Hyperplane Theorem}. It is proved in Griffiths and 
Harris \cite[p.~156]{GrHa} and Goresky and 
MacPherson~\cite[p.~153]{GoMa}.

\begin{thm} Let\/ $M$ be a compact, $m$-dimensional complex 
manifold, $N$ a nonsingular hypersurface in $M$, and\/ $L$ the 
holomorphic line bundle over $M$ associated to the divisor $N$. 
Suppose $L$ is positive. Then 
\begin{itemize}
\item[{\rm(a)}] the map $H^k(M,\C)\ra H^k(N,\C)$ induced by the 
inclusion $N\hookrightarrow M$ is an isomorphism for $0\le k\le m-2$ 
and injective for $k=m-1$, and
\item[{\rm(b)}] the map of homotopy groups $\pi_k(N)\ra\pi_k(M)$ 
induced by the inclusion $N\hookrightarrow M$ is an isomorphism for 
$0\le k\le m-2$ and surjective for~$k=m-1$.
\end{itemize}
The result also holds if\/ $M$ and\/ $N$ are orbifolds instead 
of manifolds, and\/ $N$ is a nonsingular hypersurface in the
orbifold sense.
\label{lefschetzhpthm}
\end{thm}

Here is a procedure for calculating the fundamental group and 
Betti numbers of $Y$, $Z$ and $M$. The most difficult part is
finding the Euler characteristic $\chi(Y)$, which we have already
explained above.

\begin{itemize}
\item[(a)] Calculate $\pi_1(Y)$, $H^2(Y,\C)$ and $H^3(Y,\C)$ 
explicitly. This can usually be done using Theorem 
\ref{lefschetzhpthm}. If $Y$ is a hypersurface in 
$\CP^5_{a_0,\ldots,a_5}$ then $\pi_1(Y)=\{1\}$, $H^2(Y,\C)=\C$ 
and~$H^3(Y,\C)=0$.

\item[(b)] Compute the Euler characteristic $\chi(Y)$ of $Y$,
as in~\S\ref{n63}.

\item[(c)] Calculate $\pi_1(Z)$, $H^2(Z,\C)$ and $H^3(Z,\C)$ 
from $\pi_1(Y)$, $H^2(Y,\C)$ and $H^3(Y,\C)$. Note that 
$H^j(Z,\C)$ is the $\sigma$-invariant part of $H^j(Y,\C)$.
Since $\sigma$ swaps $H^{p,q}(Y)$ and $H^{q,p}(Y)$, it follows
that $b^3(Z)=\ha b^3(Y)$.

\item[(d)] Compute the Euler characteristic $\chi(Z)$ of $Z$. 
If $\sigma$ fixes $k$ points in $Y$ then this is given 
by~$\chi(Z)=\ha\bigl(\chi(Y)+k\bigr)$.

\item[(e)] From (c) we know $b^2(Z)$ and $b^3(Z)$, and $b^1(Z)=0$ 
as $\pi_1(Z)$ is finite. Thus we can calculate $b^4(Z)$ using the 
formula~$b^4(Z)=\chi(Z)-2-2b^2(Z)+2b^3(Z)$.

\item[(f)] Now $M$ was constructed in \S\ref{n5} by gluing 
$X_{\smash{n_1}},\ldots,X_{\smash{n_k}}$ into $Z$, where
$n_j=1$ or 2 and $X_1$, $X_2$ are defined in \S\ref{n4}. 
It is easy to show that the Betti numbers of $X_1$ and $X_2$
are $b^1=b^2=b^3=0$ and $b^4=1$. Therefore the Betti numbers 
$b^j(M)$ satisfy
\begin{equation}
\!\!\!\!\!\!\!\!\!\!\!\!\!\!\!\!\!\!\!\!\!\!\!\!\!\!\!\!
b^j(M)=b^j(Z)\;\>\text{for $j=1,2,3$, and}\;\> b^4(M)=b^4(Z)+k.
\label{n6bjmeq}
\end{equation}
Also, Proposition \ref{n5prop4} gives~$\pi_1(M)$.

\item[(g)] As $M$ has metrics with holonomy $\Spin(7)$ or
$\Z_2\ltimes\SU(4)$ by Theorem \ref{n5thm2}, we know that
$\Ahat(M)=1$. Thus \eq{Ahateq} gives
\begin{equation*}
b^2(M)-b^3(M)-b^4_+(M)+2b^4_-(M)+25=0.
\end{equation*}
So we can calculate $b^4_\pm(M)$ using the equations
\begin{equation}
\begin{split}
\!\!\!\!
b^4_+(M)&={\textstyle{1\over 3}}\bigl(b^2(M)-b^3(M)+2b^4(M)+25\bigr)\\
\!\!\!\!\!\!\!\!\!\!\!\!\!\!\!\!\!\!\! \text{and}\quad
b^4_-(M)&={\textstyle{1\over 3}}\bigl(-b^2(M)+b^3(M)+b^4(M)-25\bigr).
\end{split}
\label{n6b4pmeq}
\end{equation}      
\end{itemize}

\subsection{A way of checking the answers}
\label{n65}

If you make a mistake at some stage in these calculations, which is 
quite easy to do, then you are likely not to notice unless your values 
for $b^4_\pm(M)$ are not integers. Thus it is desirable to have some 
method for checking the answers. Here is a way of doing this. All of 
our examples have been checked for consistency in this way and others, 
but for brevity we will leave out the calculations.

Suppose we can compute the Hodge number $h^{3,1}(Y)$, using complex
geometry. Then we can compute $b^4_-(Z)$ using the formula
\begin{equation*}
b^4_-(Z)=h^{3,1}(Y)+b^2(Y)-b^2(Z)-1. 
\end{equation*}
But as $X_1$ and $X_2$ have $b^4_-=1$, as in \eq{n6bjmeq} we have
$b^4_-(M)=b^4_-(Z)+k$. This gives an independent way of finding 
$b^4_-(M)$, which can be compared with your answer in part (g) above.

Now there is a complicated method for computing $h^{3,1}(Y)$
involving spectral sequences, and also a much simpler method 
called the `polynomial deformation method' which {\it does not 
always give the right answer}. Both are discussed by Green and 
H\"ubsch \cite{GrHu}. Here is a sketch of the polynomial 
deformation method.

Suppose for simplicity that $Y$ is a hypersurface of degree $d$ 
in $\CP^5_{a_0,\ldots,a_5}$. As $Y$ is a Calabi-Yau orbifold, 
$h^{3,1}(Y)$ is the dimension of the moduli space of complex 
structures on $Y$. We assume (this is {\it not} necessarily 
true) that every small deformation of $Y$ is also a hypersurface
of degree $d$ in $\CP^5_{a_0,\ldots,a_5}$, and that two nearby 
isomorphic hypersurfaces $Y,Y'$ of degree $d$ are related by an 
automorphism of~$\CP^5_{a_0,\ldots,a_5}$.

If these assumptions hold, then $h^{3,1}(Y)=m-n$, where $m$ is 
the dimension of the space of hypersurfaces of degree $d$ in 
$\CP^5_{a_0,\ldots,a_5}$, and $n$ is the dimension of the 
automorphism group of $\CP^5_{a_0,\ldots,a_5}$. Both $m$ and
$n$ are readily computed from $a_0,\ldots,a_5$ and~$d$.

\section{A simple example}
\label{n7}

Let $Y$ be the hypersurface of degree 12 in 
$\CP^5_{1,1,1,1,4,4}$ given by
\begin{equation*}
Y=\bigl\{[z_0,\ldots,z_5]\in\CP^5_{1,1,1,1,4,4}:z_0^{12}+
z_1^{12}+z_2^{12}+z_3^{12}+z_4^3+z_5^3=0\bigr\}.
\end{equation*}
Then $c_1(Y)=0$, as $12=1+1+1+1+4+4$, and $Y$ is K\"ahler as 
$\CP^5_{1,\ldots,4}$ is K\"ahler. Calculation shows that $Y$ has 
three singular points $p_1=[0,0,0,0,1,-1]$, 
$p_2=[0,0,0,0,1,e^{\pi i/3}]$ and $p_3=[0,0,0,0,1,e^{-\pi i/3}]$, 
satisfying Condition~\ref{n5cond}.

We use the method of \S\ref{n63} to calculate the Euler 
characteristic~$\chi(Y)$.

\begin{prop} The orbifold\/ $Y$ defined above has~$\chi(Y)=4887$.
\label{n7prop1}
\end{prop}

\begin{proof} Define $Y_j$ and $\pi_j$ as in \S\ref{n63}, where 
$Y_5=Y$. Then $Y_1$ is the set of 12 points $[z_0,z_1]$ in $\CP^1$ 
with $z_0^{12}+z_1^{12}=0$, and so $\chi(Y_1)=12$. Now $\pi_2:Y_2\ra
\CP^1$ is a 12-fold branched cover branched over $Y_1$, so by 
\eq{chiyjeq} we have
\begin{equation*}
\chi(Y_2)=12\chi(\CP^1)-11\chi(Y_1)=12 \cdot 2-11\cdot 12=-108.
\end{equation*}
Similarly, $\pi_3:Y_3\ra\CP^2$ is a 12-fold branched cover 
branched over $Y_2$, so that
\begin{equation*}
\chi(Y_3)=12\chi(\CP^2)-11\chi(Y_2)=12 \cdot 2-11\cdot(-108)=1224.
\end{equation*}
And $\pi_4:Y_4\ra\CP^3$ is a 3-fold branched cover of $\CP^3$ 
branched over $Y_3$, giving
\begin{equation*}
\chi(Y_4)=3\chi(\CP^3)-2\chi(Y_3)=3\cdot 4-2\cdot 1224=-2436.
\end{equation*}
Finally, $\pi_5:Y\ra\CP^4_{1,1,1,1,4}$ is a 3-fold branched cover 
of $\CP^4_{1,1,1,1,4}$ branched over $Y_4$, and so
\begin{equation*}
\chi(Y)=3\chi(\CP^4_{1,1,1,1,4})-2\chi(Y_4)=3\cdot 5-2\cdot(-2436)=4887,
\end{equation*}
as we want.
\end{proof}

\begin{prop} The Betti numbers of\/ $Y$ are
\begin{equation*}
b^0(Y)=1,\;\> b^1(Y)=0,\;\> b^2(Y)=1,\;\>
b^3(Y)=0\;\>\text{and}\;\> b^4(Y)=4883.
\end{equation*}
Also $Y\setminus\{p_1,p_2,p_3\}$ is simply-connected 
and\/~$h^{2,0}(Y)=0$.
\label{n7prop2}
\end{prop}

\begin{proof} Theorem \ref{lefschetzhpthm} shows that 
$H^k(Y,\C)\cong H^k(\CP^5_{1,\ldots,4},\C)$ for $0\le k\le 3$. Since 
$b^k(\CP^5_{1,\ldots,4})$ is 1 for $k$ even with $0\le k\le 10$ and 
0 otherwise, this shows that $b^0(Y)=b^2(Y)=1$ and $b^1(Y)=b^3(Y)=0$, 
and so $b^4(Y)=4883$ as~$\chi(Y)=4887$.

Theorem \ref{lefschetzhpthm} also gives $\pi_1(Y)\cong
\pi_1(\CP^5_{1,\ldots,4})$, so $Y$ is simply-connected. As 
the nonsingular set of $\CP^5_{1,\ldots,4}$ is simply-connected,
we can strengthen this to show that $Y\setminus\{p_1,p_2,p_3\}$ is 
simply-connected. The isomorphism $H^k(Y,\C)\cong 
H^k(\CP^5_{1,\ldots,4},\C)$ above identifies $H^{p,q}(Y)$ 
with $H^{p,q}(\CP^5_{1,\ldots,4})$, and so $h^{p,q}(Y)=
h^{p,q}(\CP^5_{1,\ldots,4})$ for $p+q\le 3$. 
Hence~$h^{2,0}(Y)=0$.
\end{proof}

Now define a map $\sigma:Y\ra Y$ by
\begin{equation*}
\sigma:[z_0,\ldots,z_5]\longmapsto[\bar z_1,-\bar z_0,\bar z_3,
-\bar z_2,\bar z_5,\bar z_4].
\end{equation*}
As in Example \ref{n6ex2}, we find that $\sigma$ is an antiholomorphic
involution of $Y$, and that the fixed points of $\sigma$ are exactly 
$p_1,p_2,p_3$. From this and \S\ref{n71} we see that Condition 
\ref{n5cond} holds for $Y$ and $\sigma$. So we can apply the 
construction of \S\ref{n5}, and resolve the orbifold $Z=Y/\an{\sigma}$ 
to get a compact 8-manifold $M$. Choosing $n_j=2$ for at least one
$j=1,2,3$, Proposition \ref{n5prop4} shows that $M$ is simply-connected,
and Theorem \ref{n5thm2} shows that $M$ admits metrics with
holonomy~$\Spin(7)$.

\begin{thm} This compact\/ $8$-manifold\/ $M$ has Betti numbers
\begin{equation*}
b^0=1,\;\> b^1=b^2=b^3=0,\;\> b^4=2446,\;\> 
b^4_+=1639\;\>\text{and}\;\> b^4_-=807.
\end{equation*}
There exist metrics with holonomy $\Spin(7)$ on $M$, which
form a smooth family of dimension~$808$.
\label{n7thm1}
\end{thm}

\begin{proof} We first calculate the Betti numbers of $Z$. 
As $\sigma$ fixes 3 points in $Y$, by properties of the Euler 
characteristic we find that $\chi(Z)=\ha(\chi(Y)+3)$. But
$\chi(Y)=4887$ by Proposition \ref{n7prop1}, so $\chi(Z)=2445$. As 
$H^k(Z,\C)$ is the $\sigma$-invariant part of $H^k(Y,\C)$ we see from 
Proposition \ref{n7prop2} that $b^0(Z)=1$ and $b^1(Z)=b^3(Z)=0$.
Also $H^2(Y,\C)$ is generated by $[\omega_{\sst Y}]$ and 
$\sigma^*(\omega_{\sst Y})=-\omega_{\sst Y}$, so $\sigma$ acts as $-1$ 
on $H^2(Y,\C)$, and~$H^2(Z,\C)=0$.

Thus $b^0(Z)=1$, $b^1(Z)=b^2(Z)=b^3(Z)=0$ and $\chi(Z)=2445$, 
giving $b^4(Z)=2443$. Equation \eq{n6bjmeq} then gives the 
Betti numbers of $M$, and \eq{n6b4pmeq} gives $b^4_\pm$. 
Theorem \ref{n5thm2} shows that there exist torsion-free 
$\Spin(7)$-structures $(\tilde\Omega,\tilde g)$ on $M$, with 
$\Hol(\tilde g)=\Spin(7)$ as $M$ is simply-connected. By Theorem 
\ref{n2thm2} the moduli space of metrics on $M$ with holonomy 
$\Spin(7)$ is a smooth manifold of dimension~$1+b^4_-(M)=808$.
\end{proof}

\subsection{A variation on this example}
\label{n71}

Here is a variation on the above, using the idea of {\it partial
crepant resolution} mentioned in \S\ref{n61}. Let $Y$ be as in 
\S\ref{n71}, but define $\sigma':Y\ra Y$ by
\begin{equation*}
\sigma':[z_0,\ldots,z_5]\longmapsto[\bar z_1,-\bar z_0,\bar z_3,
-\bar z_2,\bar z_4,\bar z_5].
\end{equation*}
Then $\sigma'$ is an antiholomorphic involution of $Y$, which fixes 
the singular point $p_1=[0,0,0,0,1,-1]$ in $Y$, and no other points. 
In particular, $\sigma'$ swaps over the other two singular 
points~$p_2,p_3$.

Thus $Y$ and $\sigma'$ do not satisfy Condition \ref{n5cond},
because the fixed set of $\sigma'$ is not the same as the singular
set $\{p_1,p_2,p_3\}$ of $Y$. To rectify this we resolve the
singular points $p_2,p_3$. Let $Y'$ be the blow-up of $Y$ 
at $p_2$ and $p_3$. Since $p_2$ and $p_3$ are modelled on
$\C^4/\an{\alpha}$, where $\alpha$ is given by \eq{n5aleq}, it
turns out that $Y'$ is a crepant resolution of $Y$, and so
is also a Calabi--Yau orbifold.

Then $Y'$ has just the one singular point $p_1$. The action
of $\sigma'$ on $Y$ lifts to $Y'$, with sole fixed point $p_1$.
Thus Condition \ref{n5cond} holds for $Y'$ and $\sigma'$.
Therefore we can apply the construction of \S\ref{n5} to 
$Y'$ and $\sigma'$, so that $Z'=Y'/\an{\sigma'}$ is a compact 
$\Spin(7)$-orbifold with one singular point $p_1$ modelled on 
$\R^8/G$. Choosing $n_1=2$ we get a resolution $M'$ of $Z'$, 
which is a compact, simply-connected 8-manifold admitting 
metrics with holonomy~$\Spin(7)$.

We shall calculate the topological invariants of $Y'$ and~$M'$.

\begin{prop} The Betti numbers of\/ $Y'$ are
\begin{equation*}
b^0=1,\;\> b^1=0,\;\> b^2=3,\;\> b^3=0 \;\>\text{and}\;\> b^4=4885,
\;\>\text{so that}\;\>\chi(Y')=4893. 
\end{equation*}
Also, $Y'\setminus\{p_1\}$ is simply-connected
and\/~$h^{2,0}(Y')=0$.
\label{n7prop3}
\end{prop}

\begin{proof} By definition $Y'$ is the blow-up of $Y$ at $p_2,p_3$. 
Each blow-up fixes $b^1$ and $b^3$ and adds 1 to $b^2$ and $b^4$. 
So the Betti numbers of $Y'$ follow from Proposition \ref{n7prop2}. 
As $Y\setminus\{p_1,p_2,p_3\}$ is simply-connected and $h^{2,0}(Y)=0$, 
we see that $Y'\setminus\{p_1\}$ is simply-connected and~$h^{2,0}(Y')=0$.
\end{proof}

Here is the analogue of Theorem~\ref{n7thm1}:

\begin{thm} This compact\/ $8$-manifold\/ $M'$ has Betti numbers
\begin{equation*}
b^0=1,\;\> b^1=0,\;\> b^2=1,\;\> b^3=0,\;\> b^4=2444,\;\>
b^4_+=1638\;\>\text{and}\;\> b^4_-=806.
\end{equation*}
There exist metrics with holonomy $\Spin(7)$ on $M'$, which
form a smooth family of dimension~$807$.
\label{n7thm2}
\end{thm}

\begin{proof} As $\sigma$ fixes 1 point in $Y'$ we have 
$\chi(Z')=\ha(\chi(Y')+1)$, so $\chi(Z')=2447$ by the previous 
proposition. Since $H^k(Z',\C)$ is the $\sigma$-invariant part of 
$H^k(Y',\C)$ we have $b^0(Z')=1$ and $b^1(Z')=b^3(Z')=0$. Now 
$b^2(Y')=3$, and $H^2(Y',\C)$ is generated by $[\omega_{\sst Y'}]$ 
and the cohomology classes Poincar\'e dual to the two exceptional 
divisors $\CP^3$ introduced by blowing up $p_2$ and $p_3$. But 
$\sigma'$ swaps $p_2$ and $p_3$, so $\sigma'_*$ swaps the corresponding 
classes in $H^2(Y',\C)$, and $\sigma'_*(\omega_{\sst Y'})=-\omega_{\sst Y'}$ 
by definition. Therefore $H^2(Y',\C)\cong\C\oplus\C^2$, where $\sigma'_*$ 
acts as 1 on $\C$ and $-1$ on $\C^2$. Hence $H^2(Z',\C)\cong\C$, 
and~$b^2(Z')=1$. 

Thus $b^0(Z')=b^2(Z')=1$, $b^1(Z')=b^3(Z')=0$ and $\chi(Z')=2447$, 
giving $b^4(Z')=2443$. Equation \eq{n6bjmeq} then gives the 
Betti numbers of $M$, and \eq{n6b4pmeq} gives $b^4_\pm$. 
Theorem \ref{n5thm2} shows that there exist torsion-free 
$\Spin(7)$-structures $(\tilde\Omega,\tilde g)$ on $M$, with 
$\Hol(\tilde g)=\Spin(7)$ as $M$ is simply-connected. By Theorem 
\ref{n2thm2} the moduli space of metrics on $M$ with holonomy 
$\Spin(7)$ is a smooth manifold of dimension~$1+b^4_-(M)=807$.
\end{proof}

Observe that the Betti numbers of $M$ and $M'$ in Theorems 
\ref{n7thm1} and \ref{n7thm2} are very similar. It is an 
interesting question whether one can regard $M$ and $M'$ as 
two different resolutions of some singular $\Spin(7)$-manifold 
$M_0$, not necessarily an orbifold. We leave this as a research 
exercise for the reader; the answer is not as simple as it looks.

\section{Examples from hypersurfaces in 
$\CP^5_{\lowercase{a_0,\ldots,a_5}}$}
\label{n8}

Here are three more examples based on hypersurfaces 
in~$\CP^5_{a_0,\ldots,a_5}$.

\subsection{A hypersurface of degree $16$ in $\CP^5_{1,1,1,1,4,8}$}
\label{n81}

Let $Y$ be the hypersurface of degree 16 in 
$\CP^5_{1,1,1,1,4,8}$ given by
\begin{equation*}
Y=\bigl\{[z_0,\ldots,z_5]\in\CP^5_{1,1,1,1,4,8}:z_0^{16}+
z_1^{16}+z_2^{16}+z_3^{16}+z_4^4+z_5^2=0\bigr\}.
\end{equation*}
Then $c_1(Y)=0$. We find that $Y$ has two singular points 
$p_1=[0,0,0,0,1,i]$ and $p_2=[0,0,0,0,1,-i]$, both satisfying
Condition~\ref{n5cond}.

Following Propositions \ref{n7prop1} and \ref{n7prop2}, we 
find that $\chi(Y)=9498$, and

\begin{prop} The Betti numbers of\/ $Y$ are
\begin{equation*}
b^0=1,\quad b^1=0,\quad b^2=1,\quad b^3=0\quad\text{and}\quad b^4=9494.
\end{equation*}
Also $Y\setminus\{p_1,p_2\}$ is simply-connected
and\/~$h^{2,0}(Y)=0$. 
\label{n8prop2}
\end{prop}

Define an antiholomorphic involution $\sigma:Y\ra Y$ by
\begin{equation*}
\sigma:[z_0,\ldots,z_5]\longmapsto[\bar z_1,-\bar z_0,\bar z_3,
-\bar z_2,\bar z_4,-\bar z_5].
\end{equation*}
The fixed points of $\sigma$ are exactly the singular points 
$p_1,p_2$ of $Y$. Thus Condition \ref{n5cond} holds for $Y$ 
and $\sigma$, and we can apply the construction of \S\ref{n5}.
Resolving $Z=Y/\an{\sigma}$ gives a compact 8-manifold $M$. We 
choose at least one of $n_1,n_2$ to be 2, so that $M$ is
simply-connected. Then as in Theorem \ref{n7thm1}, we get:

\begin{thm} This compact\/ $8$-manifold\/ $M$ has Betti numbers
\begin{equation*}
b^0=1,\;\> b^1=b^2=b^3=0,\;\> b^4=4750,\;\> 
b^4_+=3175\;\>\text{and}\;\> b^4_-=1575.
\end{equation*}
There exist metrics with holonomy $\Spin(7)$ on $M$, which
form a smooth family of dimension~$1576$.
\label{n8thm1}
\end{thm}

\subsection{A hypersurface of degree $24$ in $\CP^5_{1,1,1,1,8,12}$}
\label{n82}

Let $Y$ be the hypersurface of degree 24 in $\CP^5_{1,1,1,1,8,12}$ 
given by
\begin{equation*}
Y=\bigl\{[z_0,\ldots,z_5]\in\CP^5_{1,1,1,1,8,12}:z_0^{24}+
z_1^{24}+z_2^{24}+z_3^{24}+z_4^3+z_5^2=0\bigr\}.
\end{equation*}
Then $c_1(Y)=0$. We find that $Y$ has one singular point
$p_1=[0,0,0,0,-1,1]$, which satisfies Condition~\ref{n5cond}. 

Following Proposition \ref{n7prop1}, we find that 
$\chi(Y)=23\,325$. Care is needed to get the right answer 
here. Define $\pi_5:Y\ra\CP^4_{1,1,1,1,8}$
by $\pi_5:[z_0,\ldots,z_5]\mapsto[z_0,\ldots,z_4]$, and
$Y_4\subset\CP^4_{1,1,1,1,8}$ by
\begin{equation*}
Y_4=\bigl\{[z_0,\ldots,z_4]\in\CP^4_{1,1,1,1,8}:z_0^{24}+
z_1^{24}+z_2^{24}+z_3^{24}+z_4^3=0\bigr\}.
\end{equation*}
Then $\pi_5$ is a double cover of $\CP^4_{1,1,1,1,8}$ branched 
over $Y_4$ {\it and the point} $[0,0,0,0,1]$ in $\CP^4_{1,1,1,1,8}$.
Hence we get
\begin{equation*}
\chi(Y)=2\chi(\CP^4_{1,1,1,1,8})-\chi(Y_4)-\chi([0,0,0,0,1])
=9-\chi(Y_4).
\end{equation*}
If we had not observed that $\pi_5$ is also branched over $[0,0,0,0,1]$
then we would have got $\chi(Y)=23\,326$, which is incorrect.

As in Proposition \ref{n7prop2}, we show:

\begin{prop} The Betti numbers of\/ $Y$ are
\begin{equation*}
b^0=1,\quad b^1=0,\quad b^2=1, \quad b^3=0\quad\text{and}\quad 
b^4=23\,231.
\end{equation*}
Also $Y\setminus\{p_1\}$ is simply-connected 
and\/~$h^{2,0}(Y)=0$. 
\label{n8prop3}
\end{prop}

Define an antiholomorphic involution $\sigma:Y\ra Y$ by
\begin{equation*}
\sigma:[z_0,\ldots,z_5]\longmapsto[\bar z_1,-\bar z_0,\bar z_3,
-\bar z_2,\bar z_4,\bar z_5].
\end{equation*}
The fixed points of $\sigma$ are exactly the singular point 
$p_1$ of $Y$. Thus Condition \ref{n5cond} holds for $Y$
and $\sigma$, and choosing the simply-connected resolution 
$M$ of $Z=Y/\an{\sigma}$, in the usual way we get:

\begin{thm} This compact\/ $8$-manifold\/ $M$ has Betti numbers
\begin{equation*}
b^0=1,\;\> b^1=b^2=b^3=0,\;\> b^4=11\,662,\;\> 
b^4_+=7783\;\>\text{and}\;\> b^4_-=3879.
\end{equation*}
There exist metrics with holonomy $\Spin(7)$ on $M$, which
form a smooth family of dimension~$3880$.
\label{n8thm2}
\end{thm}

This is the example with the largest value of $b^4$ 
known to the author.

\subsection{A hypersurface of degree $40$ in $\CP^5_{1,1,5,5,8,20}$}
\label{n83}

Here is a more complicated example, in which the hypersurface
in $\CP^5_{a_0,\ldots,a_5}$ has other singularities which must
first be resolved. Let $W$ be the hypersurface of degree 40 in 
$\CP^5_{1,1,5,5,8,20}$ given by
\begin{equation*}
W=\bigl\{[z_0,\ldots,z_5]\in\CP^5_{1,1,5,5,8,20}:z_0^{40}+
z_1^{40}+z_2^8+z_3^8+z_4^5+z_5^2=0\bigr\}.
\end{equation*}
Then $c_1(W)=0$. The singularities of $W$ are the disjoint union 
of the single point $p_1=[0,0,0,0,-1,1]$ and the nonsingular curve 
$\Sigma$ of genus 3 given by
\begin{equation*}
\Sigma=\bigl\{[0,0,z_2,z_3,0,z_5]\in\CP^5_{1,1,5,5,8,20}:
z_2^8+z_3^8+z_5^2=0\bigr\}.
\end{equation*}

The singularity at $p_1$ is modelled on $\C^4/\Z_4$, where the 
generator $\alpha$ of $\Z_4$ acts on $\C^4$ by \eq{n5aleq}. The
singularity at each point of $\Sigma$ is modelled on $\C\times\C^3/\Z_5$,
where the generator $\beta$ of $\Z_5$ acts on $\C^3$ by
\begin{equation*}
\beta:(z_0,z_1,z_4)\mapsto({\rm e}^{2\pi i/5}z_0,{\rm e}^{2\pi i/5}z_1,
{\rm e}^{-4\pi i/5}z_4).
\end{equation*}
Now the singularity $\C^3/\Z_5$ normal to $\Sigma$ in $W$ has a 
unique crepant resolution $X$, which can be described using toric
geometry. Let $Y$ be the partial crepant resolution of $W$ which
resolves the singularities at $\Sigma$ using $X$, but leaves the 
singular point $p_1$ unchanged. 

\begin{prop} The Betti numbers of\/ $Y$ are
\begin{equation*}
b^0=1,\quad b^1=0,\quad b^2=3, \quad b^3=12,\quad\text{and}\quad b^4=7453.
\end{equation*}
Also $Y\setminus\{p_1\}$ is simply-connected and\/~$h^{2,0}(Y)=0$. 
\label{n8prop4}
\end{prop}

\begin{proof} Calculating the Betti numbers of $W$ in the usual 
way gives
\begin{equation}
\!\!\!\!\!\!\!\!\!\!\!\!\!\!\!\!\!\!\!\!\!\!\!\!\!\!\!\!\!\!\!
b^0(W)=1,\;\> b^1(W)=0,\;\> b^2(W)=1, \;\> b^3(W)=0,\;\> b^4(W)=7449.
\!\!\!\!\!
\label{n8wbettieq}
\end{equation}
As $W$ is modelled on $\C\times\C^3/\Z_5$ at each point of $\Sigma$, the 
resolution $Y$ is modelled on $\C\times X$. Since $b^2(X)=b^4(X)=2$, 
the Betti numbers of $Y$ satisfy
\begin{equation*}
b^k(Y)=b^k(W)+2b^{k-2}(\Sigma)+2b^{k-4}(\Sigma).
\end{equation*}
But $\Sigma$ has genus 3, and so its Betti numbers are 
$b^0(\Sigma)=b^2(\Sigma)=1$ and $b^1(\Sigma)=6$. Combining this
with \eq{n8wbettieq} gives the Betti numbers of $Y$. The last part
follows as in Proposition~\ref{n7prop2}.
\end{proof}

Define $\sigma:W\ra W$ by
\begin{equation*}
\sigma:[z_0,\ldots,z_5]\longmapsto[\bar z_1,-\bar z_0,\bar z_3,
-\bar z_2,\bar z_4,\bar z_5].
\end{equation*}
The only fixed point of $\sigma$ is $p_1$. Moreover, $\sigma$ lifts to 
the resolution $Y$ of $W$, and $\sigma:Y\ra Y$ is an antiholomorphic 
involution which fixes only $p_1$ in $Y$. Thus Condition 
\ref{n5cond} holds for $Y$ and $\sigma$, and we can apply the 
construction of \S\ref{n5}, and resolve $Z=Y/\an{\sigma}$ to get 
a simply-connected 8-manifold $M$. Proceeding in the usual way, 
the end result is

\begin{thm} This compact\/ $8$-manifold\/ $M$ has Betti numbers
\begin{equation*}
b^0=1,\;\> b^1=b^2=0,\;\> b^3=6,\;\> b^4=3730,\;\> 
b^4_+=2493\;\>\text{and}\;\> b^4_-=1237.
\end{equation*}
There exist metrics with holonomy $\Spin(7)$ on $M$, which
form a smooth family of dimension~$1238$.
\label{n8thm3}
\end{thm}

Note that $b^3>0$ in this example; this is because the
resolution of the singular curve $\Sigma$ contributes
$H^1(\Sigma,\C)\otimes H^2(X,\C)=\C^6\otimes\C^2=\C^{12}$ to 
$H^3(Y,\C)$. Half of this $\C^{12}$ is $\sigma$-invariant,
and so pushes down to $H^3(Z,\C)$ and lifts to~$H^3(M,\C)$.

\section{A hypersurface in $\CP^5_{1,1,1,1,2,2}$ over ${\mathbb Z}_2$}
\label{n9}

Let $W$ be the hypersurface of degree 8 in $\CP^5_{1,1,1,1,2,2}$ given by
\begin{equation*}
W=\bigl\{[z_0,\ldots,z_5]\in\CP^5_{1,1,1,1,2,2}:
z_0^8+z_1^8+z_2^8+z_3^8+z_4^4+z_5^4=0\bigr\}.
\end{equation*}
Then $c_1(W)=0$. We find that $W$ has four singular points 
$p_1,\ldots,p_4$ modelled on $\C^4/\{\pm1\}$, given by
\begin{equation*}
[0,0,0,0,1,{\rm e}^{\pi i/4}],\;\> [0,0,0,0,1,{\rm e}^{3\pi i/4}],\;\>
[0,0,0,0,1,{\rm e}^{5\pi i/4}],\;\> [0,0,0,0,1,{\rm e}^{7\pi i/4}].
\end{equation*}

Define $\beta:W\ra W$ by
\begin{equation*}
\beta:[z_0,\ldots,z_5]\mapsto[iz_0,iz_1,iz_2,iz_3,z_4,z_5].
\end{equation*}
Then $\beta^2=1$, as $[z_0,\ldots,z_5]=[-z_0,-z_1,-z_2,-z_3,z_4,z_5]$
in $\CP^5_{1,1,1,1,2,2}$. The fixed set of $\beta$ is the four points
$p_1,\ldots,p_4$ together with the compact complex surface $S$ in $W$,
given by
\begin{equation*}
S=\bigl\{[z_0,z_1,z_2,z_3,0,0]\in\CP^5_{1,1,1,1,2,2}:z_0^8+z_1^8+
z_2^8+z_3^8=0\bigr\}.
\end{equation*}
Thus $W/\an{\beta}$ is a compact complex orbifold. Its singular set
is the disjoint union of $p_1,\ldots,p_4$ and $S$. Each singular
point $p_j$ is modelled on $\C^4/\Z_4$, where the generator $\alpha$
of $\Z_4$ acts on $\C^4$ by \eq{n5aleq}. Each singular point in 
$S$ is locally modelled on~$\C^2\times\C^2/\{\pm1\}$.

Let $Y$ be the blow-up of $W/\an{\beta}$ along $S$. Because the
singularities normal to $S$ are modelled on $\C^2/\{\pm1\}$, this
is a {\it partial crepant resolution}. So $Y$ is a compact complex 
orbifold with isolated singular points $p_1,\ldots,p_4$, modelled on 
$\C^4/\an{\alpha}$. Now $c_1(W)=0$, so $c_1(W/\an{\beta})=0$, and as $Y$ 
is a partial crepant resolution of $W/\an{\beta}$ we see that~$c_1(Y)=0$. 

\begin{prop} The Betti numbers of\/ $Y$ are
\begin{equation*}
b^0=1,\quad b^1=0,\quad b^2=2, \quad b^3=0\quad\text{and}\quad b^4=1806.
\end{equation*}
Also $Y\setminus\{p_1,\ldots,p_4\}$ is simply-connected
and\/~$h^{2,0}(Y)=0$. 
\label{n9prop1}
\end{prop}

\begin{proof} As in Proposition \ref{n7prop1}, we find  
$\chi(W)=2708$ and $\chi(S)=304$. Thus
\begin{equation*}
\chi(W/\an{\beta})=\ha\bigl(\chi(W)+\chi(\text{4 points})+\chi(S)\bigr)
=\ha(2708+4+304)=1508.
\end{equation*}
Using Theorem \ref{lefschetzhpthm} we find that $W$ has 
$b^0=b^2=1$ and $b^1=b^3=0$, and it soon follows that
$W/\an{\beta}$ also has $b^0=b^2=1$ and $b^1=b^3=0$. Since 
$\chi(W/\an{\beta})=1508$ we see that~$b^4(W/\an{\beta})=1504$.

Now $Y$ is the blow-up of $W/\an{\beta}$ along $S$, so that each
point of $S$ is replaced by a copy of $\CP^1$. It can be shown
that the Betti numbers of $Y$ satisfy
\begin{equation}
b^k(Y)=b^k(W/\an{\beta})+b^{k-2}(S).
\label{n9ybettieq}
\end{equation}
But $S$ can be thought of as an octic in $\CP^3$, and by the
usual method we find that the Betti numbers of $S$ are
$b^0=1$, $b^1=0$, $b^2=302$, $b^3=0$ and $b^4=1$. Combining
these with \eq{n9ybettieq} and the Betti numbers of $W/\an{\beta}$
above gives the Betti numbers of $Y$. The last part follows as usual.
\end{proof}

Define an antiholomorphic involution $\sigma:W\ra W$ by
\begin{equation*}
\sigma:[z_0,\ldots,z_5]\longmapsto[\bar z_1,-\bar z_0,\bar z_3,
-\bar z_2,\bar z_5,\bar z_4].
\end{equation*}
The fixed points of $\sigma$ are exactly the singular points 
$p_1,\ldots,p_4$ of $W$. Also $\sigma$ commutes with $\beta$,
and acts freely on $S$. Thus $\sigma$ pushes down to an 
antiholomorphic involution of $W/\an{\beta}$, and lifts to 
the blow-up $Y$, to give an antiholomorphic involution
$\sigma:Y\ra Y$ with fixed points~$p_1,\ldots,p_4$.

Thus Condition \ref{n5cond} holds for $Y$ and $\sigma$, and
in the usual way we choose a simply-connected resolution
$M$ of $Z=Y/\an{\sigma}$ satisfying:

\begin{thm} This compact\/ $8$-manifold\/ $M$ has Betti numbers
\begin{equation*}
b^0=1,\;\> b^1=b^2=b^3=0,\;\> b^4=910,\;\> 
b^4_+=615\;\>\text{and}\;\> b^4_-=295.
\end{equation*}
There exist metrics with holonomy $\Spin(7)$ on $M$, which
form a smooth family of dimension~$296$.
\label{n9thm1}
\end{thm}

\subsection{A variation on this example}
\label{n91}

We shall use the idea of \S\ref{n71} to make a second 8-manifold 
$M'$ from the orbifold $Y$ above. Let $W$ and $Y$ be as in 
\S\ref{n91}, but define $\sigma':W\ra W$ by
\begin{equation*}
\sigma':[z_0,\ldots,z_5]\longmapsto[\bar z_1,-\bar z_0,\bar z_3,
-\bar z_2,\bar z_4,i\bar z_5].
\end{equation*}
Then $\sigma'$ pushes down to $W/\an{\beta}$ and lifts to $Y$ as above.
However, this time $\sigma'$ fixes the singular points $p_1=
[0,0,0,0,1,{\rm e}^{\pi i/4}]$ and $p_2=[0,0,0,0,1,{\rm e}^{5\pi i/4}]$ 
in $Y$, but it swaps round $p_3=[0,0,0,0,1,{\rm e}^{3\pi i/4}]$
and~$p_4=[0,0,0,0,1,{\rm e}^{7\pi i/4}]$.

Thus, Condition \ref{n5cond} does not hold for $Y$ and $\sigma'$, as
the fixed set $\{p_1,p_2\}$ of $\sigma'$ does not coincide with the
singular set $\{p_1,\ldots,p_4\}$ of $Y$. So let $Y'$ be the blow-up
of $Y$ at $p_3$ and $p_4$. Then $Y'$ is a partial crepant resolution 
of $Y$, as the singularities at $p_3,p_4$ are modelled on $\C^4/\Z_4$. 
The singularities of $Y'$ are $p_1,p_2$, and $\sigma'$ lifts to an 
antiholomorphic involution of $Y'$ fixing only $p_1$ and~$p_2$.

We find the Betti numbers of $Y'$ by adding contributions to 
those of $Y$, as in \S\ref{n71}. Applying the construction of 
\S\ref{n5} to $Y'$ and $\sigma'$ gives a simply-connected 8-manifold 
$M'$, such that

\begin{thm} This compact\/ $8$-manifold\/ $M'$ has Betti numbers
\begin{equation*}
b^0=1,\;\> b^1=0,\;\> b^2=1,\;\> b^3=0,\;\> b^4=908,\;\> 
b^4_+=614\;\>\text{and}\;\> b^4_-=294.
\end{equation*}
There exist metrics with holonomy $\Spin(7)$ on $M'$, which
form a smooth family of dimension~$295$.
\label{n9thm2}
\end{thm}

\section{Complete intersections in $\CP^6_{\lowercase{a_0,\ldots,a_6}}$}
\label{n10}

We now try starting with the intersection of two
hypersurfaces in~$\CP^6_{a_0,\ldots,a_6}$.

\subsection{The intersection of two octics in $\CP^6_{1,1,1,1,4,4,4}$}
\label{n101}

Let $Y$ be the complete intersection of two octics in 
$\CP^6_{1,1,1,1,4,4,4}$ given by
\begin{align*}
Y=\bigl\{[z_0,\ldots,z_5]\in\CP^6_{1,1,1,1,4,4,4}:\,
& z_0^8+z_1^8+2iz_2^8-2iz_3^8+z_4^2-z_5^2=0,\\
& 2iz_0^8-2iz_1^8+z_2^8+z_3^8+z_4^2-z_6^2=0\bigr\}.
\end{align*}
Then $c_1(Y)=0$. We find that $Y$ has 4 singular points 
$p_1=[0,0,0,0,1,1,1]$, $p_2=[0,0,0,0,1,-1,-1]$, $p_3=[0,0,0,0,1,1,-1]$
and $p_4=[0,0,0,0,1,-1,1]$ satisfying Condition~\ref{n5cond}.

By adapting the method of \S\ref{n63} we can show that
$\chi(Y)=2580$, and applying Theorem \ref{lefschetzhpthm} twice 
we find that $b^k(Y)=b^k(\CP^6_{1,\ldots,4})$ for $0\le k\le 3$.
Thus we prove:

\begin{prop} The Betti numbers of\/ $Y$ are
\begin{equation*}
b^0=1,\quad b^1=0,\quad b^2=1, \quad b^3=0\quad\text{and}\quad b^4=2576,
\end{equation*}
Also $Y\setminus\{p_1,\ldots,p_4\}$ is simply-connected and\/~$h^{2,0}(Y)=0$.
\label{n10prop1}
\end{prop}

Define an antiholomorphic involution $\sigma:Y\ra Y$ by
\begin{equation*}
\sigma:[z_0,\ldots,z_6]\longmapsto[\bar z_1,-\bar z_0,\bar z_3,
-\bar z_2,\bar z_4,\bar z_5,\bar z_6].
\end{equation*}
The fixed points of $\sigma$ are exactly the singular points 
$p_1,\ldots,p_4$ of $Y$, and Condition \ref{n5cond} holds
for $Y$ and $\sigma$. Proceeding in the usual way, we set 
$Z=Y/\an{\sigma}$ and resolve $Z$ to get a simply-connected 
8-manifold $M$, which satisfies:

\begin{thm} This compact\/ $8$-manifold\/ $M$ has Betti numbers
\begin{equation*}
b^0=1,\;\> b^1=b^2=b^3=0,\;\> b^4=1294,\;\> 
b^4_+=871\;\>\text{and}\;\> b^4_-=423.
\end{equation*}
There exist metrics with holonomy $\Spin(7)$ on $M$, which
form a smooth family of dimension~$424$.
\label{n10thm1}
\end{thm}

\subsection{A variation on this example}
\label{n102}

Now let $Y$ be as in \S\ref{n101}, but define $\sigma':Y\ra Y$ by
\begin{equation*}
\sigma':[z_0,\ldots,z_6]\longmapsto[\bar z_3,-\bar z_2,\bar z_1,
-\bar z_0,\bar z_4,\bar z_6,\bar z_5].
\end{equation*}
Then $\sigma'$ is an antiholomorphic involution, with fixed points
$p_1$ and $p_2$, which swaps round $p_3$ and $p_4$. Following the
method of \S\ref{n71}, define $Y'$ to be the blow-up of $Y$ at
$p_3$ and $p_4$. Then $Y'$ is a Calabi--Yau orbifold, $\sigma'$ lifts 
to $Y'$, and Condition \ref{n5cond} holds for $Y'$ and~$\sigma'$.

As usual we set $Z'=Y'/\an{\sigma'}$ and resolve $Z'$ to get a 
simply-connected 8-manifold $M'$, such that

\begin{thm} This compact\/ $8$-manifold\/ $M'$ has Betti numbers
\begin{equation*}
b^0=1,\;\> b^1=0,\;\> b^2=1,\;\> b^3=0,\;\> b^4=1292,\;\> 
b^4_+=870\;\>\text{and}\;\> b^4_-=422.
\end{equation*}
There exist metrics with holonomy $\Spin(7)$ on $M'$, which
form a smooth family of dimension~$423$.
\label{n10thm2}
\end{thm}

\subsection{The intersection of two 12-tics in $\CP^6_{3,3,3,3,4,4,4}$}
\label{n103}

Let $P(z_4,z_5,z_6)$ and $Q(z_4,z_5,z_6)$ be generic homogeneous
cubic polynomials with real coefficients, and define $W$ to be the 
complete intersection of two 12-tics in $\CP^6_{3,3,3,3,4,4,4}$ given by
\begin{align*}
W=\bigl\{[z_0,\ldots,z_5]\in &\,\CP^6_{3,3,3,3,4,4,4}:
z_0^4+z_1^4+z_2^4+z_3^4+P(z_4,z_5,z_6)=0,\\
&\qquad\qquad iz_0^4-iz_1^4+2iz_2^4-2iz_3^4+Q(z_4,z_5,z_6)=0\bigr\}.
\end{align*}
Then $c_1(W)=0$. As $P$ and $Q$ are generic, the singular set of $W$ 
is the disjoint union of the 9 points $p_1,\ldots,p_9$ given by
\begin{equation*}
\bigl\{[0,0,0,0,z_4,z_5,z_6]\in\CP^6_{3,3,3,3,4,4,4}:
P(z_4,z_5,z_6)=Q(z_4,z_5,z_6)=0\bigr\},
\end{equation*}
and the curve $\Sigma$ of genus 33 given by
\begin{align*}
\Sigma=\bigl\{[z_0,z_1,z_2,z_3,0,0,0]\in\CP^6_{3,3,3,3,4,4,4}:\,
& z_0^4+z_1^4+z_2^4+z_3^4=0,\\
& iz_0^4-iz_1^4+2iz_2^4-2iz_3^4=0\bigr\}.
\end{align*}

Each point $p_j$ satisfies Condition \ref{n5cond}, and each point 
of $\Sigma$ is modelled on $\C\times\C^3/\Z_3$, where the action of $\Z_3$ 
on $\C^3$ is generated by
\begin{equation*}
\beta:(z_4,z_5,z_6)\mapsto({\rm e}^{2\pi i/3}z_4,
{\rm e}^{2\pi i/3}z_5,{\rm e}^{2\pi i/3}z_6).
\end{equation*}

Define an antiholomorphic involution $\sigma:W\ra W$ by
\begin{equation*}
\sigma:[z_0,\ldots,z_6]\longmapsto[\bar z_1,-\bar z_0,\bar z_3,
-\bar z_2,\bar z_4,\bar z_5,\bar z_6].
\end{equation*}
Then the fixed points of $\sigma$ are some subset of $\{p_1,\ldots,p_9\}$.
Exactly which subset depends on the choice of $P$ and $Q$, but
$\sigma$ must fix an odd number of the $p_j$, as the remaining $p_j$
are swapped in pairs.

So let $\sigma$ fix $2k+1$ of the $p_j$, for some $k=0,\ldots,4$, and 
number the $p_j$ such that $\sigma$ fixes $p_1,\ldots,p_{2k+1}$ and 
swaps $p_{2k+2},\ldots,p_9$ in pairs. Define $Y_k$ to be the blow-up
of $W$ along $\Sigma$ and at the points $p_{2k+2},\ldots,p_9$. Then $Y_k$
is a partial crepant resolution of $W$. Thus $Y_k$ is a Calabi--Yau
orbifold, with singular points $p_1,\ldots,p_{2k+1}$. Also $\sigma$ lifts 
to $Y_k$ to give an antiholomorphic involution $\sigma:Y_k\ra Y_k$ with 
fixed points~$p_1,\ldots,p_{2k+1}$.

It can be shown that we can choose $P$ and $Q$ so that $k$ takes
any value in $\{0,1,2,3,4\}$. For example, if $P=z_4^3-z_5^3$ and 
$Q=z_4^3-z_6^3$ then $\sigma$ fixes only $p_1=[0,0,0,0,1,1,1]$, so 
that $k=0$, but if $P=z_4^2z_5-z_5^3$ and $Q=z_4^2z_6-z_6^3$
then $\sigma$ fixes the 9 points $[0,0,0,0,1,z_5,z_6]$ for
$z_5,z_6\in\{1,0,-1\}$, and~$k=4$.

Combining the methods used to prove Propositions \ref{n8prop4}
and \ref{n10prop1}, we get

\begin{prop} The Betti numbers of\/ $Y_k$ are $b^0=1$, $b^1=0$, 
$b^2=10-2k$, $b^3=66$, $b^4=395-2k$, $b^4_+=262$ and\/ $b^4_-=133-2k$.
Also $Y_k\setminus\{p_1,\ldots,p_{2k+1}\}$ is simply-connected, 
and\/~$h^{2,0}(Y_k)=0$.
\label{n10prop2}
\end{prop}

In the usual way we resolve $Z_k=Y_k/\an{\sigma}$ to get $M_k$, 
which satisfies

\begin{thm} For each\/ $k=0,\ldots,4$ there is a compact\/ 
$8$-manifold\/ $M_k$ with Betti numbers $b^0=1$, $b^1=0$, $b^2=4-k$, 
$b^3=33$, $b^4=200+2k$, $b^4_+=132+k$ and\/ $b^4_-=68+k$. There 
exist metrics with holonomy $\Spin(7)$ on $M_k$, which form a 
smooth family of dimension~$69+k$.
\label{n10thm3}
\end{thm}

These examples have the largest value of $b^3$ and the smallest 
values of $b^4$ that the author has found using this construction.

\section{Conclusions}
\label{n11}

In Table \ref{n11table} we give the Betti numbers $(b^2,b^3,b^4)$
of the compact 8-manifolds with holonomy $\Spin(7)$ that we
constructed in \S\ref{n7}--\S\ref{n10}. There are 14 sets of 
Betti numbers, none of which coincide with any in \cite{Joy1},
so we have found at least 14 topologically distinct new 
examples of compact 8-manifolds with holonomy~$\Spin(7)$.

\begin{table}[htb]
{\caption{Betti numbers $(b^2,b^3,b^4)$ of compact
$\Spin(7)$-manifolds
\qquad\qquad}\label{n11table}}
{\begin{tabular}{ccccc}
\hline
\vphantom{$\bigr)^{k^k}$}
(4,\,33,\,200) & (3,\,33,\,202) & (2,\,33,\,204) & 
(1,\,33,\,206) & (0,\,33,\,208) \\ 
(1,\,0,\,908)  & (0,\,0,\,910)  & (1,\,0,\,1292) & 
(0,\,0,\,1294) & (1,\,0,\,2444) \\ 
(0,\,0,\,2446) & (0,\,6,\,3730) & (0,\,0,\,4750) & 
(0,\,0,\,11\,662)\vphantom{$\bigr)_{p_p}$} \\
\hline
\end{tabular}}
\end{table}

The examples of \S\ref{n7}--\S\ref{n10} are by no means all 
the manifolds that can be produced using the methods of this 
paper, but only a selection chosen for their simplicity and 
to illustrate certain techniques. Readers are invited to look
for other examples themselves; the author would be particularly
interested in examples which have especially large or small
values of~$b^4$.

We have also chosen to restrict our attention in 
\S\ref{n5}--\S\ref{n10} to orbifolds $Y$ all of whose 
singularities are modelled on $\C^4/\Z_4$, where the
generator $\alpha$ of $\Z_4$ acts as in \eq{n5aleq}. This is 
not a necessary restriction, and there are other types of 
singularities for $Y$ and $Z$ for which the construction
would work, such as the $\R^8/G^n$ considered in \S\ref{n43},
and which occur in suitable orbifolds $Y$. However, the author
has not found many such $Y$; the $\C^4/\Z_4$ singularities do
seem to be the easiest to construct.

There is a very simple kind of singularity of $Y$ and $Z$ which 
can be resolved using this construction, which we have not yet
mentioned. Suppose that $Y$ has a singular point $p$ modelled 
on $\C^4/\{\pm1\}$ and fixed by $\sigma$, and that we can identify 
$T_pY$ with $\C^4/\{\pm1\}$ such that $\theta_{\sst Y}$ is identified
with $\d z_1\w\cdots\w\d z_4$ and $\d\sigma:T_pY\ra T_pY$ is identified 
with $\beta:\C^4/\{\pm1\}\ra\C^4/\{\pm1\}$, where
\begin{equation*}
\beta:(z_1,\ldots,z_4)\mapsto(\bar z_2,-\bar z_1,\bar z_4,-\bar z_3).
\end{equation*}

Then $Z=Y/\an{\sigma}$ has a singular point modelled on $\R^8/\Z_4$ 
at $p$. Changing to the complex coordinates $(w_1,\ldots,w_4)$ of 
\S\ref{n42}, the action of $\beta$ is
\begin{equation*}
\beta:(w_1,\ldots,w_4)\mapsto(iw_1,iw_2,iw_3,iw_4).
\end{equation*}
Thus the blow-up of $\C^4/\an{\beta}$ at 0 is a crepant resolution, 
as in \S\ref{n42}, and we can use this to resolve $Z$ with 
holonomy~$\Spin(7)$.

However, if the action of $\d\sigma$ on $T_pY$ is modelled
on $\beta':\C^4/\{\pm1\}\ra\C^4/\{\pm1\}$, where
\begin{equation*}
\beta':(z_1,\ldots,z_4)\mapsto(\bar z_2,-\bar z_1,-\bar z_4,\bar z_3),
\end{equation*}
then things are different. In the $w_j$ coordinates, $\beta'$ acts by
\begin{equation*}
\beta':(w_1,\ldots,w_4)\mapsto(iw_1,iw_2,-iw_3,-iw_4).
\end{equation*}
The singularity $\C^4/\an{\beta'}$ does not have a crepant resolution,
and one can in fact show that it cannot be resolved within 
holonomy~$\Spin(7)$.

Thus, the analogue of Proposition \ref{n5prop2} does not hold
for singularities modelled on $\C^4/\{\pm1\}$. When $Y$ has 
singular points $p$ modelled on $\C^4/\{\pm1\}$, there are two 
different ways $\d\sigma$ can act on $T_pY$, up to isomorphism. 
Both ways lead to $\R^8/\Z_4$ singularities in $Z$, but one can 
be resolved within holonomy $\Spin(7)$, and one cannot. 

All the Calabi--Yau orbifolds $Y$ with $\C^4/\{\pm1\}$ 
singularities that the author has looked at, such as the octic 
in $\CP^5_{1,1,1,1,2,2}$, seem to have at least one singular point 
of each kind, and so cannot be resolved with this construction.

\end{document}